\documentclass{amsart}
\usepackage{amsmath,amssymb,amscd,latexsym}
\usepackage[all]{xy}

\newcommand{\C}{{\mathbb{C}}}

\newcommand{\N}{{\mathbb{N}}}

\newcommand{\R}{{\mathbb{R}}}

\newcommand{\Ch}{{\mathcal C}}

\newcommand{\Kh}{{\mathcal K}}

\newcommand{\Ph}{{\mathcal P}}

\newcommand{\Zh}{{\mathcal Z}}

\newcommand{\Aff}{\mathrm{Aff}}

\newcommand{\be}{\mathbf{1}}

\newcommand{\cpr}{\mathrm{cpr }\,}

\newcommand{\dist}{\mathrm{dist}}
\newcommand{\dr}{\mathrm{dr}\,}

\newcommand{\halb}{\frac{1}{2}}
\newcommand{\her}{\mathrm{her}}

\newcommand{\id}{\mathrm{id}}
\newcommand{\ord}{\mathrm{ord}\,}

\newcommand{\rr}{\mathrm{rr}\,}

\newcommand{\verk}{\mbox{\scriptsize $\,\circ\,$}}

\newcounter{number}[subsection]
\newcounter{altnumber}[section]

\newenvironment{nummer}{\refstepcounter{number}{\noindent\arabic{section}.\arabic{subsection}.\arabic{number}}}{}
\newenvironment{altnummer}{\refstepcounter{altnumber}{\noindent\arabic{section}.\arabic{altnumber}}}{}

\newcommand{\bn}{\noindent\begin{nummer} \rm}
\newcommand{\en}{\end{nummer}}

\newcommand{\altbn}{\noindent \begin{altnummer} \rm}
\newcommand{\alten}{\end{altnummer}}

\newenvironment{ntheorem}{\noindent {\sc Theorem:} \it}{}
\newenvironment{nlemma}{\noindent {\sc Lemma:} \it}{}
\newenvironment{nprop}{\noindent {\sc Proposition:} \it}{}
\newenvironment{ndefn}{\noindent {\sc Definition:} \it}{}
\newenvironment{ncor}{\noindent {\sc Corollary:} \it}{}
\newenvironment{nconj}{\noindent {\sc Conjecture:} \it}{}
\newenvironment{nremark}{\noindent {\sc Remark:} }{}
\newenvironment{nremarks}{\noindent {\sc Remarks: }}{}
\newenvironment{nexamples}{\noindent {\sc Examples:} }{}

\newenvironment{nproof}{\noindent {\sc Proof:}}{\mbox{}\hfill \rule[-.2ex]{.25em}{1.8ex}}

\begin{document}

\title[Topologically finite-dimensional simple $C^{*}$-algebras]{{\sc On topologically finite-dimensional \\
simple $C^{*}$-algebras}}

\author{Wilhelm Winter}
\address{Mathematisches Institut der Universit\"at M\"unster\\ 
Einsteinstr. 62\\ D-48149 M\"unster}

\email{wwinter@math.uni-muenster.de}

\date{November 2003}
\thanks{{\it Supported by:} EU-Network Quantum Spaces - Noncommutative Geometry (Contract No. HPRN-CT-2002-00280) and Deutsche Forschungsgemeinschaft (SFB 478)}

\setcounter{section}{-1}

\begin{abstract}
We show that, if a simple $C^{*}$-algebra $A$ is topologically finite-dimensional in a suitable sense, then not only $K_{0}(A)$ has certain good properties, but $A$ is even accessible to Elliott's classification program. More precisely, we prove the following results: \\
If $A$ is simple, separable and unital with finite decomposition rank and real rank zero, then $K_{0}(A)$ is weakly unperforated. \\
If $A$ has finite decomposition rank, real rank zero and the space of extremal tracial states is compact and zero-dimensional, then $A$ has stable rank one and tracial rank zero. As a consequence, if $B$ is another such algebra, and if $A$ and $B$ have isomorphic Elliott invariants and satisfy the Universal coefficient theorem, then they are isomorphic. \\
In the case where $A$ has finite decomposition rank and the space of extremal tracial states is compact and zero-dimensional, we also give a criterion (in terms of the ordered $K_{0}$-group) for $A$ to have real rank zero. As a byproduct, we show that there are examples of simple, stably finite and quasidiagonal $C^{*}$-algebras with infinite decomposition rank. 
\end{abstract}

\maketitle

\section{Introduction}
\noindent
The theory of $C^{*}$-algebras is often regarded as noncommutative topology. This point of view  is already suggested by Gelfand's Theorem; however, it has experienced  great impact since the introduction of methods from (algebraic) topology to $C^{*}$-algebra theory. In particular, (bivariant) $K$-theory has turned out to be a strong tool to study $C^{*}$-algebras. In the present article we combine various notions of noncommutative covering dimension with $K$-theoretic methods to obtain classification results for certain simple nuclear $C^{*}$-algebras.  \\
In \cite{El1},  Elliott showed that $K$-theory provides a complete invariant for a certain natural class of $C^{*}$-algebras, the $AF$-algebras.  This classification result has become a model for what is now known as the Elliott program (see \cite{Ro} for an introduction); it is the aim of this program to find complete $K$-theoretic invariants, such as $K$-theory, the trace simplex and the natural pairing between these objects, for suitable classes of nuclear $C^{*}$-algebras. \\
In these notes, we shall only deal with simple and stably finite $C^{*}$-algebras (a $C^{*}$-algebra $A$ is simple if it has no nontrivial closed ideals; $A$ is stably finite if none of the matrix algebras $M_{k}(A)$ contains a projection which is infinite in the sense of Murray--von Neumann).  Moreover, for most of our results we will restrict ourselves to unital algebras with real rank zero (i.e., self-adjoint elements with finite spectrum are dense in the set of all self-adjoint elements).  For a nuclear, simple, stably finite, unital $C^{*}$-algebra $A$ the classifying invariant suggested by Elliott is the tuple
\[
(K_{0}(A),K_{0}(A)_{+},[\be_{A}],T(A),r_{A}:T(A) \to S(K_{0}(A)),K_{1}(A)) \, .
\]
Here, $(K_{0}(A),K_{0}(A)_{+},[\be_{A}])$ is the ordered $K_{0}$-group with distinguished order unit $[\be_{A}]$ and state space $S(K_{0}(A))$; $T(A)$ is the space of tracial states on $A$ and $r_{A}:T(A) \to S(K_{0}(A))$ the natural map given by the pairing between $T(A)$ and $K_{0}(A)$. For nuclear $A$ the map $r_{A}$ is surjective; if $A$ has real rank zero, $r_{A}$ is also injective (cf.\ \cite{Ro}, 1.1.11 and 1.1.12). For such  $A$, the invariant is only  $(K_{0}(A),K_{0}(A)_{+},[\be_{A}],K_{1}(A))$ and the Elliott conjecture reads as follows:

\altbn{\label{elliott-conjecture}}
\begin{nconj}
Let $\mathfrak{A}$  be the class of separable, simple, unital, stably finite and nuclear $C^{*}$-algebras with real rank zero and let $A$, $B \in \mathfrak{A}$. Then $A$ and $B$ are isomorphic if and only if they have isomorphic Elliott invariants. Moreover, every isomorphism of the invariants is induced by an isomorphism of the algebras.
\end{nconj}
\alten

\noindent
We say \ref{elliott-conjecture} holds for a subclass $\mathfrak{B}$ of $\mathfrak{A}$ if the assertion holds for all $A, \, B \in \mathfrak{B}$. The subclasses  of $\mathfrak{A}$ for which \ref{elliott-conjecture} has been verified so far consist of  approximately homogeneous ($AH$)  algebras. These are inductive limits of $C^{*}$-algebras built from matrix valued functions over topological spaces.  Using the covering dimension of these underlying spaces, one can easily define a notion of topological dimension for each inductive limit decomposition (and then also for the limit algebra). $AH$-algebras with finite topological dimension in this sense are particularly accessible to classification, cf.\ \cite{Ro}. One can also define the topological dimension of approximately subhomogeneous ($ASH$) algebras along these lines. However, for $AH$-algebras (which may also be regarded as $ASH$-algebras) these two definitions do not coincide in general.

\noindent
In the absence of an obvious inductive limit decomposition it is less clear how to define the topological dimension of a $C^{*}$-algebra. In \cite{KW}, E.\ Kirchberg and the author introduced the decomposition rank of a nuclear $C^{*}$-algebra $A$, $\dr A$. This is defined by imposing a certain condition on systems of completely positive (c.p.) approximations for $A$ (see Section 1 for details). The decomposition rank  generalizes  topological covering dimension to nuclear $C^{*}$-algebras and has good permanence properties; it is computable for many examples (cf.\ \cite{KW} and \cite{Wi3}).
 
\noindent
There are various other notions of covering dimension for $C^{*}$-algebras, such as stable and real rank (see \cite{Ri1} and \cite{BP}, respectively), both of which have proven to be extremely important for the classification program.\\
Villadsen was the first to construct simple $AH$-algebras with higher stable and real ranks; these examples have recently been modified by Toms who showed that there are simple $AH$-algebras (with nonzero real rank) for which the Elliott invariant is not complete (the invariant does not detect the stable rank of the algebras). These  examples have infinite  topological dimension (as $AH$-algebras); together with the available classification results for $AH$-algebras they suggest that the Elliott conjecture is much better accessible to $C^{*}$-algebras which are topologically finite-dimensional in a suitable sense. 

\noindent
To begin with, we will show that a simple unital $C^{*}$-algebra $A$ which has real rank zero and  finite decomposition rank  satisfies Blackadar's second fundamental comparability property - in other words, projections  in $A$ may be compared in terms of tracial states on $A$. This is possible, since one can approximate the elements  and, at the same time, the tracial states of $A$  by a suitably chosen  system of completely positive (c.p.) approximations. As an immediate consequence, the ordered group $K_{0}(A)$ is weakly unperforated (cf.\ \ref{d-w-u}). \\
To obtain a classification theorem, we will not verify \ref{elliott-conjecture} directly; instead, we use the concept  of $C^{*}$-algebras with tracial  rank zero. Roughly speaking, $A$ has tracial rank zero if it contains a  sequence of finite-dimensional subalgebras that `almost' (i.e., in terms of traces) exhaust $A$. This notion was introduced by Lin to approach Conjecture \ref{elliott-conjecture} for algebras which do not have an obvious  inductive limit decomposition. In \cite{Li3}, Lin has shown that tracial rank zero algebras   indeed satisfy \ref{elliott-conjecture} -- provided the Universal coefficient theorem holds (cf.\ \cite{Bl1}). \\
In order to verify that our algebras have tracial rank zero,  we have to ask for an extra condition -- similar to real rank zero -- on the tracial state space $T(A)$. The technical reason is, that we have to be able to approximate certain partitions of unity of $T(A)$ by pairwise orthogonal elements of the algebra. This is possible, if the extreme boundary $\partial_{e}T(A)$ of the compact convex set $T(A)$ is compact and zero-dimensional. Our main theorem then is the following:

\altbn
\label{main-theorem}
\begin{ntheorem} Let $A$ be a separable, simple and unital $C^{*}$-algebra with finite decomposition rank and real rank zero. Suppose that $\partial_{e} T(A)$ is compact and zero-dimensional. Then $A$ has tracial rank zero.
\end{ntheorem}           
\alten

\noindent 
From this we can deduce a number of Corollaries. First note that it follows from Lin's work that Conjecture \ref{elliott-conjecture} holds for algebras as in Theorem \ref{main-theorem}, if they  also satisfy the Universal coefficient theorem. Such algebras are $AH$ of topological dimension at most 3; they are $ASH$ with decomposition rank (and topological dimension) at most 2.\\
For algebras as in the theorem it also turns out that they are invariant under tensoring with the Jiang--Su algebra $\Zh$. This is an algebra $KK$-equivalent to $\C$; for algebras as in \ref{elliott-conjecture} which have weakly unperforated $K_{0}$-groups, the Elliott conjecture predicts that these are tensorially invariant w.r.t.\ $\Zh$. \\
If $A$ is simple and unital with finite decomposition rank and if $\partial_{e}T(A)$ is compact and zero-dimensional, we show that $A$ has real rank zero if and only if $K_{0}(A)_{+}$ has dense image in  $\Ch(\partial_{e}T(A))$ under the natural map. If $A$ has only one tracial state $\tau$, then $A$ has real rank zero if and only if $\tau(K_{0}(A))$ is dense in $\R$; otherwise, $A$ is stably isomorphic to a unital $C^{*}$-algebra with no nontrivial projections. We then apply  this dichotomy to  show that certain  examples of Villadsen (which are simple $AH$ with infinite topological dimension and perforated $K_{0}$-groups)  have infinite decomposition rank; this answers a question left open in \cite{KW}.\\
Finally, we show that indeed there are many $C^{*}$-algebras satisfying the conditions of Theorem \ref{main-theorem}.

\noindent
This paper benefitted greatly from discussions with N.\ Brown, J.\ Cuntz, S.\ Eilers, N.\ C.\ Phillips, M.\ R{\o}rdam and, in particular, E.\ Kirchberg. Part of this work was done while the author was a visiting assistant professor at Texas A\&M University; the author is indebted to the members of that institution and especially to Ken Dykema for their kind hospitality.

\newpage

\section{Decomposition rank and maps of strict order zero}

\noindent
Below we recall the definitions of decomposition rank and of strict order zero maps; we generalize the partial ordering on the positive cone of a $C^{*}$-algebra to such maps. Furthermore, we prove a technical result about simple $C^{*}$-algebras with finite decomposition rank and outline a relation between approximating systems of a $C^{*}$-algebra and its tracial state space.
 
\altbn{\label{d-dr}}
Recall that nuclear $C^*$-algebras are characterized by the completely positve approximation property, i.e., $A$ is nuclear if and only if there is a net of finite-dimensional $C^*$-algebras $F_\lambda$ and completely positive contractive (c.p.c.) maps $A \stackrel{\psi_\lambda}{\longrightarrow} F_\lambda \stackrel{\varphi_\lambda}{\longrightarrow} A$ such that $\varphi_\lambda \verk \psi_\lambda$ converges to $\id_A$ pointwise (see \cite{Tak} for an introduction to nuclear $C^{*}$algebras and completely positive maps). We then say $(F_\lambda,\psi_\lambda,\varphi_\lambda)_\Lambda$ is a system of c.p.\ approximations for $A$. Based on this approximation property, E.\ Kirchberg and the author have defined a noncommutative version of covering dimension as follows:

\begin{ndefn}
(cf.\ \cite{KW}, Definitions 2.2 and 3.1) Let $A$ be a separable $C^*$-algebra. \\
(i) A c.p.\ map $\varphi : \bigoplus_{i=1}^s M_{r_i} \to A$ has strict order zero, $\ord \varphi = 0$, if it preserves orthogonality, i.e., $\varphi(e) \varphi(f) = \varphi(f) \varphi(e) = 0$ for all $e,f \in \bigoplus_{i=1}^s M_{r_i}$ with $ef = fe = 0$.\\ 
(ii) A c.p.\ map $\varphi : \bigoplus_{i=1}^s M_{r_i} \to A$ is $n$-decomposable, if there is a decomposition $\{1, \ldots, s\} = \coprod_{j=0}^n I_j$ s.t.\ the restriction of $\varphi$ to $\bigoplus_{i \in I_j} M_{r_i}$ has strict order zero for each $j \in \{0, \ldots, n\}$; we say $\varphi$ is $n$-decomposable w.r.t.\ $\coprod_{j=0}^n I_j$.\\
(iii) $A$ has decomposition rank $n$, $\dr A = n$, if $n$ is the least integer such that the following holds: Given $\{b_1, \ldots, b_m\} \subset A$ and $\varepsilon > 0$, there is a c.p.\ approximation $(F, \psi, \varphi)$ for $b_1, \ldots, b_m$ within $\varepsilon$ (i.e., $\psi:A \to F$ and $\varphi:F \to A$ are c.p.c.\ and $\|\varphi \psi (b_i) - b_i\| < \varepsilon$) such that $\varphi$ is $n$-decomposable. If no such $n$ exists, we write $\dr A = \infty$.  
\end{ndefn}

\noindent
This notion has good permanence properties; for example, it behaves well with respect to quotients, inductive limits, hereditary subalgebras, unitization and stabilization. Furthermore, it generalizes topological covering dimension, i.e., if $X$ is a locally compact second countable space, then $\cpr \Ch_0(X) = \dr \Ch_0(X) = \dim X$; see \cite{KW} for details. 
\alten

\altbn{\label{order-zero}}
In \cite{Wi1}, Proposition 4.4.1(a), maps of strict order zero were characterized as follows: \\
If $\varphi : F \to A$ is c.p.c.\ with $\ord \varphi = 0$, then there is a unique $*$-homomorphism $\pi_\varphi : CF \to A$ such that $\pi_\varphi(g \otimes x) = \varphi(x) \; \forall \, x \in F$, where $CF$ is the cone $\Ch_0((0,1]) \otimes F$ over $F$ and $g := \id_{(0,1]}$ is the canonical generator of $\Ch_0((0,1])$. Conversely, any $*$-homomorphism $\pi : CF \to A$ induces such a c.p.c.\ map $\varphi$ of strict order zero.\\
The $*$-homomorphism $\pi_\varphi$ extends to a $*$-homomorphism $\pi_\varphi'' : (CF)'' \to (\varphi(F))'' \subset A''$ of von Neumann algebras. Denote by $\sigma : F \to A''$ the $*$-homomorphism coming from the composition of the natural unital embedding $F \hookrightarrow (CF)''$ and $\pi_\varphi''$. Then we have $\varphi (x) = \varphi(\be_F) \sigma (x) = \sigma(x) \varphi(\be_F) \; \forall \, x \in F$; we call any $*$-homomorphism : $ F \to A''$ with this property a supporting $*$-homomorphism for $\varphi$. 

\noindent
Note that each $h \in C^*(\varphi(\be_F))' \cap A$ satisfying $0 \le h \le \be$ and $h \sigma(F) \subset A$, defines a c.p.c.\ map $\hat{\varphi}: F \to A$ by $\hat{\varphi}(\,.\,):= h \sigma(\, .\,)$ s.t.\ $\ord \hat{\varphi} = 0$ and $\| \hat{\varphi} - \varphi\| = \|h - \varphi(\be_F)\|$. 
\alten

\altbn{\label{induced-order0}}
Let $\sigma: M_{r} \to A''$ be a $*$-homomorphism and $e\in M_{r}$ a minimal projection. Suppose $d \in A$ is a positive element such that $d \le \sigma(e)$ (thus $\sigma(e)d \in A$). Then $d$ and $\sigma$ induce a c.p.c.\ map $\varphi:M_{r} \to A$ by the formula $\varphi(e_{i,j})=\sigma(e_{i,1})d\sigma(e_{1,j})$, where $\{e_{i,j} \, | \, 1 \le i,j \le r \}$ is a set of matrix units for $M_{r}$ with $e_{1,1}=e$ (it is straightforward to check that $\varphi$ indeed maps $M_{r}$ to $A$). $\varphi$ has strict order zero,  supporting $*$-homomorphism $\sigma$ and satisfies $\varphi(e)=d$. Moreover, if $\bar{\varphi}:M_{r}\to A$ is another c.p.c.\ strict order zero map with supporting $*$-homomorphism $\sigma$ and $\bar{\varphi}(e)=d$, then $\bar{\varphi}=\varphi$. \\
In particular, if $\varphi$ is as above and $0\le d' \in \overline{dAd}$ with $\|d'\|\le 1$, there is a unique strict order zero map $\varphi':M_{r } \to A$ with supporting $*$-homomorphism $\sigma$ and $\varphi'(e)=d'$.\\
Note that $\varphi'(x)=\varphi'(\be_{M_{r}}) \sigma(x) \; \forall \, x \in M_{r}$; it is then straightforward to check that $\varphi'(x) \in \overline{\varphi(x)A\varphi(x)} \; \forall \, x \in (M_{r})_{+}$. Furthermore, if $d' \le d$, then $\varphi'(x)\le \varphi(x) \; \forall \, x \in (M_{r})_{+}$. 
\alten

\altbn{\label{subordinate}}
\begin{ndefn}
Let $\varphi: F=M_{r_{1}}\oplus \ldots \oplus M_{r_{s}} \to A$ be c.p.\ with $\ord(\varphi|_{M_{r_{i}}})=0 \; \forall \, i$. If $\varrho:F \to A$ is another c.p.\ map, we say $\varrho$ is subordinate to $\varphi$, if $\varrho(x) \in \overline{\varphi(x)A\varphi(x)} \; \forall \, x \in F_{+}$ (this implies that $\ord(\varrho|_{M_{r_{i}}})=0 \; \forall \, i$) and if, for each $i$,  the restrictions $\varrho|_{M_{r_{i}}}$ and $\varphi|_{M_{r_{i}}}$ have common supporting $*$-homomorphisms. 
\end{ndefn}

\begin{nremarks}
(i) Let $\varphi:F \to A$ be as in the definition; denote the restrictions $\varphi|_{M_{r_{i}}}$ by $\varphi_{i}$ and the respective supporting $*$-homomorphisms by $\sigma_{i}$. For $i=1, \ldots, s$, let $e_{i}\in M_{r_{i}}$ be rank-one projections and let $d_{i} \in \overline{\varphi_{i}(e_{i})A\varphi(e_{i})}$ be positive normalized elements. We then clearly have $d_{i} \le \sigma_{i}(e_{i})$ and $d_{i} \sigma_{i}(e_{i})=d_{i}\in A$, so by \ref{induced-order0} we obtain induced order-zero maps $\varrho_{i}:M_{r_{i}} \to A$. The $\varrho_{i}$ now add up to a (not necessarily contractive) c.p.\ map $\varrho:F \to A$ which obviously is subordinate to $\varphi$.\\
(ii) If $\varphi:F \to A$ is $n$-decomposable and $\varrho$ is subordinate to $\varphi$, then $\varrho$ is also $n$-decomposable.\\
(iii) It is clear from the definition that, if $\varrho$ is subordinate to $\varphi$ and $\varrho'$ is subordinate to $\varrho$, then $\varrho'$ is subordinate to $\varphi$.
\end{nremarks}
\alten

\altbn{\label{weakly-stable}}
Recall from \cite{Wi2}, 1.2.3, that order zero maps are weakly stable:\\
\begin{nprop} For $F=M_{r_{1}} \oplus \ldots \oplus M_{r_{s}}$ and $\gamma>0$ there is $\alpha>0$ such that the following holds:\\
If $A$ is a $C^{*}$-algebra and $\varphi: F \to A$ a c.p.c.\ map with $\ord \varphi|_{M_{r_{i}}} = 0 \; \forall \, i$ and $\|\varphi(\be_{M_{r_{i}}}) \varphi(\be_{M_{r_{j}}})\|< \alpha$ whenever $i \neq j$, then there is a c.p.c.\ map $\bar{\varphi}:F \to A$ with $\ord \bar{\varphi} = 0$ and $\|\varphi - \bar{\varphi}\|< \gamma$.
\end{nprop}
\alten

\altbn{\label{nf-approximations}}
It was shown in \cite{KW}, Proposition 5.2, that if $A$ is a separable $C^{*}$-algebra with $\dr A=n<\infty$,  there is a system $(F_{\nu},\psi_{\nu},\varphi_{\nu})_{\nu \in \N}$ of c.p.\ approximations for $A$ such that the  $\varphi_{\nu}$ are $n$-decomposable and the $\psi_{\nu}$ are approximately multiplicative at the same time. In particular, we not only have $\|\varphi_{\nu}\psi_{\nu}(x)-x\| \stackrel{\nu \to \infty}{\longrightarrow}0$, but also $\|\psi_{\nu}(xy)-\psi_{\nu}(x)\psi_{\nu}(y)\| \stackrel{\nu \to \infty}{\longrightarrow}0$ for all $x$ and $y$ in $A$. In \cite{KW} this observation was used to show that finite decompostion rank implies quasidiagonality. 
\alten

\altbn{\label{simple-nf-approximations}} {\label{functions}} The next proposition shows that, if $\dr A< \infty$ and if $A$ happens to be simple, the approximations may be chosen with arbitrarily large matrix blocks. Recall that a $C^{*}$-algebra is nonelementary if it is not isomorphic to  the compact operators on some Hilbert space.\\
We shall use the following notation from \cite{Wi2}, 3.2.4: For positive numbers $\alpha\le \beta$  we define  continuous positive functions on $\R$ by
\[
f_{\alpha,\beta} (t) :=  \left\{ 
    \begin{array}{cl}
0 & \mbox{for} \; t \le \alpha \\
t & \mbox{for} \; \beta  \le t \\
\mbox{linear} &  \mbox{elsewhere} 
\end{array} \right.
\]
and
\[
g_{\alpha,\beta} (t) :=  \left\{ 
    \begin{array}{cl}
0 & \mbox{for} \; t \le \alpha \\
1 & \mbox{for} \; \beta  \le t \\
\mbox{linear} &  \mbox{elsewhere} \, .
\end{array} \right.
\]
Let $\chi_{\alpha}$ denote the characteristic function of $(\alpha , \infty)$.

\begin{nprop}
Let $A$ be simple and nonelementary with $\dr A=n<\infty$. Let $d\in A$ be a positive normalized element and $0\le \alpha < 1$. \\
Then there is a system of c.p.\ approximations $(F_{k},\psi_{k},\varphi_{k})_{\N}$ for $A$ such that the $\psi_{k}$ are approximately multiplicative, the $\varphi_{k}$ are $n$-decomposable and, for each $k$, the irreducible representations of $h_{k}F_{k}h_{k}$, where $h_{k}:=\chi_{\alpha}(\psi_{k}(d))$, are at least $k$-dimensional.
\end{nprop}

\begin{nproof}
It clearly suffices to prove that, given $k\in \N$, there are a system $(F_{l}, \psi_{l},\varphi_{l})_{\N}$ of c.p.\ approximations for $A$ and  $N \in \N$ such that the irreducible representations of $h_{l}F_{l}h_{l}$ are at least $k$-dimensional for $l\ge N$.\\
We start with an arbitrary system of c.p.\ approximations $(F'_{l},\psi'_{l},\varphi'_{l})_{\N}$ for $A$ with approximately multiplicative  $\psi'_{l}$ and $n$-decomposable $\varphi'_{l}$.\\
Set $d':= f_{\alpha,1}(d) \neq 0$. By \cite{Wi2}, 1.2.3 in connection with \cite{Wi2}, Lemma 3.2.2, for any $k \in \N$ there is a nonzero map $\varrho: M_{k} \to \overline{d'Ad'}$ which has  strict order zero. We may clearly assume that  $\varrho(\be_{M_{k}})\le d'$.\\
Since the $\psi'_{l}$ are approximately multiplicative, they induce a $*$-homomorphism $\psi: A \to \prod F'_{l}/ \bigoplus F'_{l}$; $\psi$ is injective because $\varphi'_{l} \verk \psi'_{l} \to \id_{A}$ pointwise. \\
By \cite{Wi2}, Proposition 1.2.4, the order zero map $\psi \verk \varrho$ has a (nonzero) c.p.c.\ lift $\bar{\psi}: M_{k} \to \prod F'_{l}$ of strict order zero. For each $l$, we may write $F'_{l}$ in the form $E_{l} \oplus F_{l}$, where $F_{l}$ consists of those direct summands of $F'_{l}$ for which the respective summand of $\bar{\psi}$ is nonzero. We have $\prod F'_{l}= \prod E_{l} \oplus \prod F_{l}$; for each $l$, let $\bar{\psi}_{E_{l}}$ and $\bar{\psi}_{F_{l}}$ denote the respective summands of $\bar{\psi}$.  Write $\psi$ as $\psi_{E} \oplus \psi_{F}$ in the obvious way, then  $\bar{\psi}_{E_{l}}=0 \; \forall \, l$ implies $\psi_{E} \verk \varrho(\be_{M_{k}}) =0$. But $\varrho(\be_{M_{k}})\neq 0$, so by simplicity of $A$ we see that $\psi_{E}=0$, whence $(F_{l},\psi_{l}:=\be_{F_{l}}\psi'_{l}(\, . \, ), \varphi_{l}:=\varphi'_{l}|_{F_{l}})_{\N}$ is a system of c.p.\ approximations for $A$ with approximately multiplicative $\psi_{l}$ and $n$-decomposable $\varphi_{l}$. Finally, we have 
\[
h_{l}=  \chi_{\alpha}(\psi_{l}(d)) \ge f_{\alpha,1}(\psi_{l}(d))\, , 
\]
\[
 \| f_{\alpha,1}( \psi_{l}(d)) -  \psi_{l}(d')\| \to 0
\]
and 
\[
\|\bar{\psi}_{l}(\be_{M_{k}})- \psi_{l}(\varrho(\be_{M_{k}}))\| \to 0 \, .
\]
Since   $\varrho(\be_{M_{k}})\le d'$ and $\bar{\psi}_{l}(\be_{M_{k}})F_{l}\bar{\psi}_{l}(\be_{M_{k}})$ has rank at least $k$ for all $l$, it follows  that $h_{l}F_{l}h_{l}$ has rank at least $k$ if only $l$ is large enough.
\end{nproof}
\alten

\altbn{\label{order0-traces}}
Suppose $\tau \in T(A)$ is a tracial state on $A$ and $\varphi: F \to A$ is an order-zero map, then $\tau \verk \varphi$ is a trace on $F$. This is because, for any $x \in F$, 
\begin{eqnarray*}
\tau \verk \varphi(x^{*}x) & = & \tau(\varphi(\be_{F}) \sigma(x^{*}) \sigma(x)) \\
& = & \tau(\varphi(\be_{F}) \sigma(x) \sigma(x)^{*}) \\
& = & \tau \verk \varphi (x x^{*}) \, ,
\end{eqnarray*}
where $\sigma$ is a supporting $*$-homomorphism for $\varphi$. The statement still holds if $\varphi$ is $n$-decomposable, i.e., a sum of order-zero maps.
\alten

\altbn{\label{qd-traces}}
If $\dr A$ is finite and $A$ is unital, then $A$ is quasidiagonal by \cite{KW} and we know from   \cite{Vo}, that $A$ has a tracial state. Applying Voiculescu's idea to our setting, we obtain a  description of tracial states in terms of approximating systems for $A$:\\
For each $\nu \in \N$ fix a tracial state $\tau_{\nu}$ on $F_{\nu}$; note that $\tau_{\nu}$ is a convex combination of the (unique) tracial states on the matrix blocks of $F_{\nu}$. Choose a free ultrafilter $\omega$ on $\N$ and define a linear map $\tau:A \to \C$ by
\[
\tau(x):= \lim_{\nu \to \omega} \ \tau_{\nu} \verk \psi_{\nu}(x) \; \forall \, x\in A\,.
\]
It is clear that $\tau$ is positive, since all the $\psi_{\nu}$ and $\tau_{\nu}$ are. Moreover, if $A$ is unital, then $\tau(\be_{A})=1$, so $\tau$ is a state on $A$.  Finally, we check that
\begin{eqnarray*}
|\tau(x^{*}x)-\tau(xx^{*})| &=& \lim_{\nu \to \omega} |\tau_{\nu}(\psi_{\nu}(x^{*}x)-\psi_{\nu}(xx^{*}))|  \\
&=&\lim_{\nu \to \omega} |\tau_{\nu}(\psi_{\nu}(x^{*})\psi_{\nu}(x)-\psi_{\nu}(x)\psi(x^{*}))| \\
&=& 0\, ,
\end{eqnarray*}
where for the second equality it is essential that $\omega$ is a {\it free} ultrafilter so we can apply approximate multiplicativity of the $\psi_{\nu}$. For the last equality we used that the $\tau_{\nu}$ are traces, i.e.\ $\tau_{\nu}(y^{*}y)=\tau_{\nu}(yy^{*})$ for any $y \in F_{\nu}$. But this means that $\tau$ actually is a tracial state on $A$. 

\noindent
We should also mention that in fact every tracial state on $A$ arises in the manner described above and that the tracial state space $T(A)$ can be written as an inverse limit of (geometric realizations of) finite simplical complexes associated to the approximating system $(F_{\nu},\psi_{\nu},\varphi_{\nu})_{\N}$.  Similar techniques then show that the tracial states on $A$ satisfy the approximation properties of \cite{Br1}. This will be made precise in subsequent work.
\alten

\section{Real rank zero}

\noindent
In this section we introduce discrete order zero maps to construct special approximating systems of $C^{*}$-algebras with finite decomposition rank and real rank zero.

\altbn{\label{rr0}}
Recall from \cite{BP} that a $C^{*}$-algebra $A$ is said to have real rank zero, $\rr A =0$, if the set of self-adjoint elements with finite spectrum is dense within the set of all self-adjoint elements in $A$. This in particular implies that, given $a \in A_{+}$ and $\varepsilon >0$, there is a positive element with finite spectrum $h\in \overline{aAa}$ such that $\|a - h\|<\varepsilon$, where $\overline{aAa} \subset_{\her} A$ denotes the hereditary $C^{*}$-subalgebra of $A$ generated by $a$. Replacing $h$ by $\tilde{h}:= (h -\varepsilon\cdot\be)_{+} \in \overline{aAa}$, we obtain an element with finite spectrum such that $\|a- \tilde{h}\|<2 \varepsilon$ and $0 \le \tilde{h}\le a$.  (By \cite{BP}, Corolary 2.8, any hereditary $C^{*}$-subalgebra of $A$ has real rank zero if $A$ has.)
\alten

\altbn{\label{discretemaps}}
\begin{ndefn}
(i) Let $A$ be a $C^{*}$-algebra. We say a c.p.\ map $\varphi:  F=M_{r_{1}}\oplus \ldots \oplus M_{r_{s}} \to A$ is a discrete order zero map, if $\ord \varphi =0$ and each $\varphi(\be_{M_{r_{i}}})$ is a multiple of a projection. \\
(ii) A c.p.\ map $\varphi:F \to A$ is discretely $n$-decomposable, if  $F$ can be written as $F=F^{(0)}\oplus \ldots \oplus F^{(n)}$ with $\varphi|_{F^{(i)}}$ being a discrete order zero map for all $i=0,\ldots,n$. 
\end{ndefn}
\alten

\altbn{\label{discretesigma}}
Let $\varphi$ be a discrete order zero map as above. Define a central element $d:=\sum_{i}\|\varphi(\be_{M_{r_{i}}})\|\cdot \be_{M_{r_{i}}} \in F_{+}$ and a c.p.\ map $\sigma_{\varphi}:F \to A$ by $\sigma_{\varphi}(x):=\varphi(d^{-1}x), \, x\in F$, where the inverse is taken in $dFd \subset F$. Note that $\sigma_{\varphi}(d)=\varphi(\be_{F})$. Since $d$ is central, $\sigma_{\varphi}$ has strict order zero; by construction, $\sigma_{\varphi}(\be_{F}) \in A$ is a projection. But then it follows from \cite{Wi1}, Proposition 3.2, that $\sigma_{\varphi}$ is a $*$-homomorphism - in fact, $\sigma_{\varphi}$ is a supporting $*$-homomorphism for $\varphi$. Clearly, 
\[
\varphi(x)=\sigma_{\varphi}(dx)=\varphi(\be_{F}) \sigma_{\varphi}(x)\; \forall \, x\in F
\]
and $[\varphi(\be_{F}),\varphi(F)]=[\varphi(\be_{F}),\sigma_{\varphi}(F)]=0$.
\alten

\altbn{\label{rr0discrete}}
\begin{nlemma}
Let $A$ and $F'$ be $C^{*}$-algebras, $A$ with real rank zero and $F'$ finite-dimensional. Suppose $\varphi':F' \to A$ is c.p.c.\ with strict order zero and let $\eta>0$ be given. Then there are a unital embedding $\iota:F' \to F$ of $F'$ into some finite-dimensional $C^{*}$-algebra $F$ and a discrete order zero map $\varphi:F \to A$ such that $\varphi(\be_{F})\le \varphi'(\be_{F'})$ and $\|\varphi'(x) - \varphi \verk \iota(x)\|< \eta \cdot \|x\|$ for all $x \in F'$. 
\end{nlemma}

\begin{nproof}
Let $F'=M_{r_{1}}\oplus \ldots \oplus M_{r_{s}}$ be given; set  $\varphi'_{i}:=\varphi'|_{M_{r_{i}}}$, $i=1,\ldots,s$. Define  $A_{i}:=\overline{\varphi'_{i}(\be_{M_{r_{i}}})A\varphi'_{i}(\be_{M_{r_{i}}})} \subset_{\her} A$; the $A_{i}$ have real rank zero since they are hereditary subalgebras of $A$ (cf.\ \ref{rr0}).  Furthermore, $A_{i}\perp A_{j}$ if $i\neq j$ and the maps  $\varphi'_{i}:M_{r_{i}} \to A_{i}$ are c.p.c. with strict order zero.  Suppose the assertion of the lemma holds in case $F'$ is a single matrix algebra. This means we can construct finite-dimensional $C^{*}$-algebras $F_{i}$, unital embeddings $\iota_{i}:M_{r_{i}} \to F_{i}$ and discrete order zero maps $\varphi_{i}:F_{i} \to A_{i}$ with $\varphi_{i}(\be_{F_{i}})\le \varphi'_{i}(\be_{M_{r_{i}}})$ and $\|\varphi'_{i}(x)-\varphi_{i} \verk \iota_{i} (x)\|<\eta \cdot \|x\| \; \forall \, x \in M_{r_{i}}$. Then put $F:=\bigoplus F_{i}$, $\iota:=\bigoplus \iota_{i}:F' \to F$ and $\varphi:= \bigoplus_{i} \varphi_{i}: F \to \bigoplus A_{i} \subset A$. But now (using that the $A_{i}$ are mutually orthogonal) one checks that \begin{eqnarray*}
\|\varphi'(x) - \varphi \verk \iota(x)\| & = & \|\sum \varphi'_{i}(\be_{M_{r_{i}}}x) - \sum \varphi_{i} \verk \iota_{i} (\be_{M_{r_{i}}}x)\| \\
& \le & \max_{i} \{ \| \varphi'_{i}(\be_{M_{r_{i}}}x) -  \varphi_{i} \verk \iota_{i} (\be_{M_{r_{i}}}x)\| \} \\
& \le & \max_{i} \{ \eta \cdot \|\be_{M_{r_{i}}}x\| \} \\
& = & \eta \cdot \|x\| 
\end{eqnarray*} 
for all $x \in F'$; it is clear that $\varphi(\be_{F})\le \varphi'(\be_{F'})$. Therefore, we only have to prove the lemma for $F'=M_{r}$. 

\noindent
So let $r\in \N$ be given. Note that, if $\varphi':M_{r} \to A$ has strict order zero, then $\varphi'(x)=\varphi'(\be_{M_{r}})\sigma(x)$, where $\sigma:M_{r}\to A''$ is the $*$-homomorphism from \ref{order-zero}. \\
Take a set $\{e_{k,l}\,|\,k,l=1,\ldots, r\}$ of matrix units for $M_{r}$ and define $a:=\varphi'(e_{11})$. Now by \ref{rr0} there are pairwise orthogonal projections $q_{1},\ldots,q_{m}\in \overline{aAa}$ and $0<\lambda_{1},\ldots,\lambda_{m}\le1$  such that  $\|\sum_{i} \lambda_{i} \cdot q_{i} - a\|<\eta$; we may also assume that $\sum_{i} \lambda_{i} \cdot q_{i} \le a$. Note that $\sigma(e_{1,1})q_{i}=q_{i} \; \forall \, i$, because $\sigma(e_{1,1})a=a$ and the $q_{i}$ live in $\overline{aAa}$. Moreover, it is not hard to check that the $q_{i}$ can be written as $q_{i}=\varphi'(\be_{M_{r}}) h_{i} \varphi'(\be_{M_{r}})$ for suitable elemets $h_{i}\in A$.

\noindent
For $i=1,\ldots,m$, define
\[
p_{i}:=\sum_{j=1}^{r} \sigma(e_{j,1})q_{i}\sigma(e_{1,j}) \,.
\]
Because $\sigma$ is a $*$-homomorphism and since $\sigma(e_{1,1})q_{i}=q_{i}$, we see that the $p_{i}$ are mutually orthogonal projections. Since $q_{i}=\varphi'(\be_{M_{r}}) h_{i} \varphi'(\be_{M_{r}})$ and $\varphi'(x)=\sigma(x)\varphi'(\be_{M_{r}})\in A$ for all $x\in M_{r}$, the $p_{i}$ in fact live in $A$. Moreover, $[p_{i},\sigma(e_{k,l})]=0 \; \forall \, k, l =1,\ldots,r$, so by \ref{order-zero} we can define strict order zero maps  $\varphi_{i}:M_{r}\to p_{i}Ap_{i}$ by $\varphi_{i}(x):= \lambda_{i} \cdot p_{i} \sigma(x)$ for $x \in M_{r}$. \\
Define $F:=\bigoplus_{i=1}^{m} M_{r}$ and $\varphi:F \to \bigoplus_{i} p_{i}Ap_{i}$ by $\varphi:=\bigoplus_{i} \varphi_{i}$. Since the $\varphi_{i}$ have pairwise orthogonal images and each $\varphi_{i}$ has strict order zero, we see that $\varphi$ also has strict order zero. Moreover,  $\varphi_{i}(\be_{M_{r}})= \lambda_{i} \cdot p_{i}$, so in fact $\varphi$ is a discrete order zero map. Note that $\varphi(\be_{F})=\sum_{i} \lambda_{i} \cdot p_{i}$, so
\begin{eqnarray*}
\|\varphi(\be_{F}) - \varphi'(\be_{M_{r}})\| &=& \|\sum_{i}\lambda_{i} \cdot p_{i}- \varphi'(\be_{M_{r}})\| \\
& = &\|\sum_{i=1}^{m} \lambda_{i} (\sum_{j=1}^{r} \sigma(e_{j,1})q_{i}\sigma(e_{1,j})) -  \sum_{j=1}^{r} \sigma(e_{j,1}) \varphi'(e_{1,1}) \sigma(e_{1,j})\| \\
& = & \| \sum_{j=1}^{r} \sigma(e_{j,1}) (\sum_{i=1}^{m} \lambda_{i}\cdot q_{i}) \sigma(e_{1,j}) - \sum_{j=1}^{r}  \sigma(e_{j,1}) a \sigma(e_{1,j})\| \\ 
& = & \| \sum_{j=1}^{r} \sigma(e_{j,1})  (\sum_{i=1}^{m} \lambda_{i}\cdot q_{i} -a) \sigma(e_{1,j})\| \\
&=& \max_{j} \|\sigma(e_{j,1})  (\sum_{i=1}^{m} \lambda_{i}\cdot q_{i} -a) \sigma(e_{1,j})\| \\
&<& \eta \, .
\end{eqnarray*}
Similarly, one checks that $\varphi(\be_{F})\le \varphi'(\be_{M_{r}})$. Let $\iota:M_{r}\to F$ be the diagonal embedding, then, for any $x \in M_{r}$, $\varphi \verk \iota(x)= \varphi(\be_{F}) \sigma(x)$, hence
 \[
 \| \varphi \verk \iota(x) - \varphi'(x)\| = \| (\varphi(\be_{F})-\varphi'(\be_{M_{r}})) \sigma(x)\| < \eta \cdot \|x\| \, .
 \]
\end{nproof}
\alten

\altbn{\label{rr0dr}}
\begin{nprop}
Let $A$ be a $C^{*}$-algebra with real rank zero and decomposition rank $n$. For any $b_{1},\ldots,b_{k} \in A$ and $\varepsilon>0$, there is a c.p.\ approximation $(F,\psi,\varphi)$ such that $\varphi$ is discretely $n$-decomposable and 
\[
\|\varphi\psi(b_{j})-b_{j}\|,\, \|\psi(b_{j})\psi(b_{l})-\psi({b_{j}b_{l}})\|<\varepsilon  \mbox{ for } j,\, l=1,\ldots,k \, .
\]
\end{nprop}

\begin{nproof}
We may assume that the $b_{j}$ are positive and normalized. Choose a c.p.\ approximation $(F',\psi',\varphi')$ such that $\varphi'$ is $n$-decomposable and
\[
\|\varphi'\psi'(b_{j})-b_{j}\|,\, \|\psi'(b_{j})\psi'(b_{l}) \psi'({b_{j}b_{l}})\|<\frac{\varepsilon}{2} \mbox{ for } j,\, l=1,\ldots,k \, ;
\]
this is possible by \cite{KW}, Proposition 5.2. Since $\varphi'$ is $n$-decomposable, we can write $F'=F'_{0}\oplus\ldots\oplus F'_{n}$ such that the restrictions $\varphi'|_{F'_{i}}$ have strict order zero for $i=0,\ldots,n$. Set $\eta:=\frac{\varepsilon}{2(n+1)}$ and apply Lemma \ref{rr0discrete} to each of these maps. This yields finite-dimensional $C^{*}$-algebras $F_{i}$, unital embeddings $\iota_{i}:F'_{i}\to F_{i}$ and discrete order zero maps $\varphi_{i}:F_{i} \to A$ (for $i=0,\ldots,n$) satisfying 
\[
\|\varphi'(x) - \varphi_{i} \verk \iota_{i}(x)< \eta \|x\| \mbox{ for all } x \in F'_{i} \, . 
\]
By \ref{rr0discrete} we may even assume that $\varphi_{i}(\be_{F_{i}})\le \varphi'(\be_{F'_{i}})$. 

\noindent
Now set $F:=F_{0}\oplus \ldots \oplus F_{n}$, $\psi:= (\iota_{0} \oplus \ldots \oplus \iota_{n}) \verk \psi'$ and $\varphi:=\varphi_{0} + \ldots + \varphi_{n}$,  then $\varphi$ is  discretely $n$-decomposable by our construction. Note that not only $\psi$, but also $\varphi$ is contractive since $\varphi(\be_{F}) \le \sum \varphi'(\be_{F'_{i}}) = \varphi'(\be_{F'})$ and $\varphi'$ is contractive. Furthermore, 
\[
\|\psi(b_{j})\psi(b_{l})-\psi({b_{j}b_{l}})\| \le \|\psi'(b_{j})\psi'(b_{l}) \psi'({b_{j}b_{l}})\| <\varepsilon
\]
because $\iota_{0} \oplus \ldots \oplus \iota_{n}$ is a $*$-homomorphism. Finally,
\begin{eqnarray*}
\|\varphi \psi(b_{j}) -\varphi'\psi'(b_{j})\|&=&\| \sum_{i} \varphi_{i} \verk \iota_{i} (\be_{F'_{i}}\psi'(b_{j}))  - \varphi'(\be_{F'_{i}}\psi'(b_{j}))\| \\
&\le&\sum_{i} \eta \|\be_{F'_{i}}\psi'(b_{j})\| \\
&\le&(n+1)\eta \, ,
\end{eqnarray*}
so
\[
\|\varphi\psi(b_{j})-b _{j}\| < \frac{\varepsilon}{2} +\frac{\varepsilon}{2} =\varepsilon \, .
\]
\end{nproof}
\alten

\altbn{\label{simple-rr0dr}}
\begin{nremark}
In the preceding proof, if the matrix blocks of $F'$ at least have size $r$, then so have the matrix blocks of $F$. Therefore, in the setting of Proposition \ref{simple-nf-approximations}, we may additionally assume the approximations $(F_{k},\psi_{k},\varphi_{k})$ to be discretely $n$-decomposable, if $A$ has real rank zero.
\end{nremark}
\alten

\section{Comparison}

\noindent
In this section we first recall some well-known facts about the $K_{0}$-group and the tracial state space of a $C^{*}$-algebra. We then prove that a simple $C^{*}$-algebra with finite decomposition rank and real rank zero has Blackadar's second fundamental comparability property, in other words, it admits comparison of projections in terms of tracial states.

\altbn
Recall that, for a unital $C^{*}$-algebra $A$, the abelian group $K_{0}(A)$ is given by
\[
K_{0}(A) = \{[p]_{0}-[q]_{0} \, | \, p,q \in \Ph(A \otimes \Kh) \} \, .
\]
Here, $\Ph(A \otimes \Kh)$ denotes the set of projections in $\Kh \otimes A$ and $[ \,.\,]_{0}$ an equivalence relation on $\Ph(A \otimes \Kh)$ based on Murray--von Neumann equivalence (see \cite{Bl1} for a detailed definition). With the positive cone 
\[
K_{0}(A)_{+} = \{[p]_{0} \, | \, p \in \Ph(A \otimes \Kh) \}  \, ,
\]
$(K_{0}(A),K_{0}(A)_{+})$ becomes a preordered group in the sense of \cite{Go}. We will only deal with stably finite $A$, in which case $(K_{0}(A),K_{0}(A)_{+})$ is an ordered  group, i.e.\, $K_{0}(A)_{+} \cap -K_{0}(A)_{+}= \emptyset$. Note that the class $[\be_{A}]$ is a distinguished order unit of $K_{0}(A)$. From now on, we will simply write $K_{0}(A)$ for the triple $(K_{0}(A),K_{0}(A)_{+},[\be_{A}])$.
\alten

\altbn{\label{r_A}}
By $S(K_{0}(A))$ we denote the state space of $K_{0}(A)$, i.e., the  set of all positive homomorphisms $f: K_{0}(A) \to \R$ with $f([\be_{A}])=1$ equipped with the weak*-topology. With this topology, $S(K_{0}(A))$ is a compact convex set, just as the space $T(A)$ of tracial states on $A$. There is a natural map $r_{A}: T(A) \to S(K_{0}(A))$. We are mainly interested in nuclear $C^{*}$-algebras of real rank zero; for these the map $r_{A}$ is a homeomorphism (cf.\ \cite{Ro}, 1.1.11 and 1.1.12). To simplify notation, we will not always write the $r_{A}$ explicitly.\\
By evaluation, the map $r_{A}$ induces a map $\iota: K_{0}(A) \to \Aff (T(A))$, the space of continuous affine functions on the convex set $T(A)$. Sometimes we shall consider the composition of $\iota$ with the restriction map $\Aff(T(A)) \to \Ch(\partial_{e}T(A))$, where $\partial_{e}T(A)$ denotes the extreme boundary of $T(A)$. This composition will also be denoted by $\iota$.
\alten

\altbn{\label{comparability}}
The space $S(K_{0}(A))$ can be used to compare elements of $K_{0}(A)$, but sometimes the space $T(A)$ (and, therefore, $S(K_{0}(A))$) can be used to compare elements (or at least projections) of $A$. More precisely: A $C^{*}$-algebra $A$ is said to have Blackadar's second fundamental comparability property, if, for any $p,q \in \Ph(A)$, the estimate $\tau(p) < \tau(q) \; \forall \, \tau \in T(A)$ implies $p \preceq q$ (i.e., $p$ is Murray--von Neumann equivalent to a subprojection of $q$). In this case, 
\[
K_{0}(A)_{+}= \{ x \in K_{0}(A) \, | \, r_{A}(\tau)(x)>0 \; \forall \, \tau \in T(A) \} \, .
\]
\alten 

\altbn{\label{d-w-u}}
The ordered group $K_{0}(A)$ is said to be weakly unperforated if it is unperforated for the strict order. This means, whenever $g\in K_{0}(A)$, $n \in \N$ and $0 \neq n \cdot g \in K_{0}(A)_{+}$, then $g \in K_{0}(A)_{+}$. 
\alten

\altbn{\label{comparability-w-u}}
\begin{nprop}
Suppose $A$ is unital and simple and $M_{k}(A)$ satisfies Blackadar's second fundamental comparability property for each $k\in \N$; let $T(A)$ be nonempty. Then $K_{0}(A)$ is weakly unperforated. 
\end{nprop}

\begin{nproof}
Suppose $g\in K_{0}(A)$, $n \in \N$ and $0 \neq n \cdot g \in K_{0}(A)_{+}$. Then there are  nonzero projections $e,f,p \in M_{k}(A)$ such that $[e]_{0}-[f]_{0}=g$ and $[p]_{0} = n \dot g =n \dot [e]_{0}-n \cdot [f]_{0}$. Now, for any $\tau \in T(A)$, since $A$ is simple we have $0< \tau(p)$, whence $0<n(\tau(e)-\tau(f))$. But then $\tau(f)<\tau(e)$ for all $\tau$ and comparison yields $f \preceq e$. This in turn means that there is $q \in M_{k}(A)$ with $[q]_{0}=[e]_{0}-[f]_{0}$, so $g \in K_{0}(A)_{+}$ and $K_{0}(A)$ is weakly unperforated.  
\end{nproof}
\alten

\noindent
In order to prove the main result of this section (Theorem \ref{perforation}) we first need a number of technical intermediate results.

\altbn{\label{interpolation}}
\begin{nlemma}
Let $A$ be a simple unital $C^{*}$-algebra with $\rr A=0$ and $\dr A=n< \infty$. if $p\in A$ is a nonzero projection, there is a projection $q\in pAp$ such that $\frac{1}{9(n+1)^{2}} \cdot \tau(p)<\tau(q)<\halb \cdot \tau(p)$ for all tracial states $\tau \in T(A)$. 
\end{nlemma}

\begin{nproof}
Set $B:=pAp$ and note that $\rr B=0$ and $\dr B=n$ by \cite{BP}, Theorem 2.5 and \cite{KW}, Proposition 3.10. Choose $\eta>0$ such that $8 (n+1) \eta<\frac{1}{9}$. \\
By Remark \ref{simple-rr0dr}, there is a c.p.\ approximation $(F=M_{r_{1}}\oplus \ldots \oplus M_{r_{s}}, \psi,\varphi)$ for $B$ with $\|\varphi \psi(p)-p\|<\eta$,  $\varphi$ discretely $n$-decomposable and such that $r_{i}>8(n+1)$ for $i=1,\ldots,s$. Note that $\|\varphi \psi(p)-p\|<\eta$ and the fact that $\varphi$ and $\psi$ are contractive together imply that $\|\varphi(\be_{F})-p\|<\eta$. \\
For each $i$, choose $l_{i}\in \N$ such that $\frac{1}{8(n+1)}<\frac{l_{i}}{r_{i}}<\frac{1}{4(n+1)}$ and a projection $e_{i} \in M_{r_{i}}$ with rank $l_{i}$ (such $l_{i}$ and $e_{i}$ obviously exist).  \\
Define elements $g_{i} \in B$ by 
\[
g_{i}:=\frac{1}{\|\varphi(\be_{M_{r_{i}}})\|} \cdot \varphi(e_{i}) 
\]
(w.l.o.g.\ we may assume that $\varphi|_{M_{r_{i}}}\neq 0 \; \forall \, i$). Since $\varphi$ is discretely $n$-decomposable, we see that the $g_{i}$ are projections. Moreover, there is a decomposition of the index set $\{1,\ldots,s\}$ into $n+1$ pairwise disjoint subsets $I_{j}$, $\{1,\ldots,s\}= \coprod_{j=0}^{n} I_{j}$, such that the restriction of $\varphi$ to $\bigoplus_{i\in I_{j}} M_{r_{i}}$ has strict order zero for each $j$. But this means that $g_{i} \perp g_{i'}$ if $i, i'\in I_{j}$ for some $j$, from which in turn follows that 
\[
q_{j}:= \sum_{i\in I_{j}} g_{i}
\]
are projections in $B$ for $j=0, \ldots,n$.\\
Next, choose some $0<\varepsilon<\frac{1}{4}$. Since $\rr B=0$, there is a projection $q \in B$ with 
\[
\halb \cdot q \le (1-2 \varepsilon) \cdot q \le g:= \sum_{j=0}^{n} q_{j}
\]
and 
\[
q\ge f_{1-\varepsilon,1}(g)\, ,
\]
where $f_{1-\varepsilon,1}$ is defined as in \ref{functions}. Note that $\|g-f_{1-\varepsilon,1}(g)\|\le 1-\varepsilon$, so, for $j=0,\ldots,n$, 
\[
q_{j}qq_{j}\ge q_{j} f_{1-\varepsilon,1}(g)q_{j} \ge q_{j}gq_{j} - (1-\varepsilon)\cdot q_{j} \ge q_{j}-(1-\varepsilon) \cdot q_{j} = \varepsilon \cdot q_{j}\, .
\]
As a consequence, $q_{j}qq_{j}$ is invertible in $q_{j}Bq_{j}$ and 
\[
q_{j}\preceq q \; \forall \, j \, . 
\]
Now let $\tau$ be some tracial state on $A$ (then $0<\tau(p)\le1$). For each $i=1,\ldots,s$, $\varphi(M_{r_{i}})$ is isomorphic to $M_{r_{i}}$ and the restriction of $\tau$ to $\varphi(M_{r_{i}})$ is a trace by \ref{order0-traces}. Therefore, 
\[
\tau(g_{i})= \frac{l_{i}}{r_{i}} \tau(\be_{M_{r_{i}}})
\]
for all $i$. It follows that 
\[ 
\tau(q_{j})> \frac{1}{8(n+1)} \cdot \tau \left( \sum_{i\in I_{j}} \frac{1}{\|\varphi(\be_{M_{r_{i}}})\|} \varphi(\be_{M_{r_{i}}}) \right) > \frac{1}{8(n+1)} \cdot \tau \left( \sum_{i\in I_{j}}\varphi(\be_{M_{r_{i}}}) \right) \, .
\]
Using that $q_{j}\preceq q$, we obtain 
\begin{eqnarray*} 
(n+1) \tau(q) & \ge & \sum_{j=0}^{n} \tau(q_{j}) \\
&>& \frac{1}{8(n+1)} \cdot \tau \left( \sum_{j=0}^{n} \sum_{i\in I_{j}}\varphi(\be_{M_{r_{i}}}) \right) \\
&=& \frac{1}{8(n+1)} \tau(\varphi(\be_{F}))\\
&\ge& \frac{1}{8(n+1)} (1-\eta) \tau(p) \, ,
\end{eqnarray*}
whence
\[
\tau(q) > \frac{1-8 (n+1)\eta}{8(n+1)^{2}} \tau(p) > \frac{1}{9(n+1)^{2}} \tau(p)\, .
\]
Conversely,
\begin{eqnarray*} 
\tau(q) &\le& 2 \tau(g)\\
& = & 2 \sum_{j=0}^{n} \tau(q_{j})\\
&<& 2 \frac{1}{4(n+1)} \sum_{j=0}^{n} \sum_{i\in I_{j}} \tau \left( \frac{1}{\|\varphi(\be_{M_{r_{i}}})\|} \varphi(\be_{M_{r_{i}}}) \right) \\
&\le&   \frac{1}{2(n+1)} \sum_{j=0}^{n} \tau(p) \\
& \le &  \frac{1}{2(n+1)} (n+1) \tau(p)\\
&=&\halb \tau(p) \, .
\end{eqnarray*}
\end{nproof}
\alten

\altbn{\label{orthogonal-subprojections}}
\begin{ncor}
Given $n\in \N$ and $\eta>0$ there are $k\in \N$ (with $2^{k-1}\ge n+1$) and $\gamma>0$ such that the following holds:\\
Let  $A$ be a simple unital $C^{*}$-algebra with $\rr A=0$ and $\dr A=n$ and let $p\in A$ be a nonzero projection. \\
Then there are mutually orthogonal projections $q_{1},\ldots ,q_{2^{k}} \in pAp$  with $p =q_{1}\oplus \ldots \oplus q_{2^{k}}$ and $\gamma \cdot \tau(p) < \tau(q_{i})<\eta \cdot \tau(p)$ for all $\tau \in T(A)$ and $i=1,\ldots,2^{k}$.  
\end{ncor}

\begin{nproof}
Choose $k \in \N$ such that $2^{k-1}\ge n+1$ and $(1-\frac{1}{9(n+1)^{2}})^{k}<\eta$. Set $\gamma:= (\frac{1}{9(n+1)^{2}})^{k}$.\\
Now let $A$ and $p$ be as in the hypothesis of the Corollary. By Lemma \ref{interpolation}, there is a projection $q \le p$ in $A$ such that, for all tracial states $\tau$ on $A$, 
\[
\frac{1}{9(n+1)^{2}} \tau(p)<\tau(q)<\halb \tau(p) \, .
\]
Note that these inequalities also imply
\[
\halb \tau(p)<\tau(p-q)< \left( 1-\frac{1}{9(n+1)^{2}} \right) \tau(p)
\]
for all $\tau \in T(A)$. Next, apply Lemma \ref{interpolation} to $q$ and to $p-q$ to obtain projections $q'\le q$ and $q''\le p-q$. Proceed inductively in the same manner; after $k$ steps we have pairwise orthogonal projections $q_{1}, \ldots, q_{2^{k}} \in pAp$ satisfying $q_{1}\oplus \ldots \oplus q_{2^{k}} =p$ and (using our inductive construction and the fact that $\frac{1}{9(n+1)^{2}}<\halb<1-\frac{1}{9(n+1)^{2}}$)
\[
\gamma \cdot \tau(p) = \left( \frac{1}{9(n+1)^{2}}\right)^{k} \tau(p) < \tau(q_{i})< \left( 1-\frac{1}{9(n+1)^{2}} \right)^{k} \tau(p) <\eta \cdot \tau(p)
\]
for $i=1, \ldots, 2^{k}$ and all $\tau \in T(A)$.  
\end{nproof}
\alten
   
\altbn{\label{n+1-domination}}
\begin{nprop}
Let $A$ be a separable, simple and unital $C^{*}$-algebra with $\dr A=n<\infty$. Suppose $p\in A$ is a projection and $d^{(0)}, \ldots, d^{(n)} \in A$ are pairwise orthogonal positive normalized elements in $A$ satisfying $\tau(p)< \tau(d^{(j)})$ $\forall \, \tau \in T(A), \, j=0, \ldots,n$. \\
Then $p$ is Murray--von Neumann equivalent to a projection in $\overline{dAd}$, where $d=d^{(0)}+\ldots +d^{(n)}$.
\end{nprop}
\alten

\begin{nproof}
Choose an approximating system $(F_{k},\psi_{k},\varphi_{k})_{\N}$ for $A$ such that the $\psi_{k}$ are approximately multiplicative and the $\varphi_{k}$ are $n$-decomposable with respect to decompositions $F_{k}=\bigoplus_{j=0}^{n}F_{k}^{(j)}$ and $\varphi_{k}^{(j)}:F_{k}^{(j)}\to A$.\\
Define 
\[\textstyle
B:= \prod_{\N}A/\bigoplus_{\N}A\, ,\; C:= \prod_{\N}F_{k}/\bigoplus_{\N}F_{k}  \mbox{ and } C^{j}:= \prod_{k \in \N} F_{k}^{(j)}/\bigoplus_{k \in \N} F_{k}^{(j)} 
\]
for $j=0, \ldots,n$; note that $C= \bigoplus_{j=0}^{n}C^{(j)}$. Let $\iota:A \to B$ be the natural inclusion (as constant sequences). The $\psi_{k}$ and $\varphi_{k}$ induce c.p.c.\ maps $\bar{\psi}: B \to C$ and $\bar{\varphi}:C \to B$. Since the $\psi_{k}$ are approximately multiplicative, $\bar{\psi}|_{\iota(A)}$ is a $*$-homomorphism; furthermore, $\bar{\varphi} \verk \bar{\psi} \verk \iota= \iota$, since $\varphi_{k}\verk \psi_{k}\to \id_{A}$ pointwise. The maps $\bar{\psi}$ and $\bar{\varphi}$ can be decomposed into maps $\bar{\psi}^{(j)}:B \to C^{(j)}$ and $\bar{\varphi}^{(j)}:C^{(j)} \to B$; clearly, the $\bar{\psi}^{(j)}$ are $*$-homomorphisms as well. \\
Since each $\ord \varphi_{k}^{(j)}=0$ for each $j$ and $k$, there are supporting $*$-homomorphisms $\sigma_{k}^{(j)}:F_{k}^{(j)} \to A''$ (cf.\ \ref{order-zero}) such that
\[
\varphi_{k}^{(j)}(x)=\varphi_{k}^{(j)}(\be_{F_{k}^{(j)}}) \sigma_{k}^{(j)}(x)=\sigma_{k}^{(j)}(x) \varphi_{k}^{(j)}(\be_{F_{k}^{(j)}})\; \forall\, x \in F_{k}^{(j)}\, .
\] 
These induce $*$-homomorphisms $\bar{\sigma}^{(j)}: C^{(j)}\to \prod_{\N} A''/\bigoplus_{\N}A''$ with the property that 
\[
\bar{\varphi}^{(j)}(x)= \bar{\varphi}^{(j)}(\be_{C^{(j)}})\bar{\sigma}^{(j)}(x)=\bar{\sigma}^{(j)}(x) \bar{\varphi}^{(j)}(\be_{C^{(j)}}) \in B\, ;
\]
here, we have identified $B$ with its image in $\prod A''/\bigoplus A''$. 

\noindent
Having introduced all this notation, we are now prepared to give the actual proof of the proposition. For each $j$, the function $T(A) \ni \tau \mapsto \tau(d_{j})-\tau(p)$ is nonzero and continuous on $T(A)$; since $T(A)$ is compact, the number 
\[
\alpha:= \min\{\tau(d_{j})-\tau(p) \, | \, j=0, \ldots, n, \, \tau \in T(A)\}
\]
exists and is nonzero. \\
For each $k\in \N$ and $j=0,\ldots, n$ consider the projections $q_{k}^{(j)}:= \chi_{\frac{\alpha}{2}}(\psi_{k}^{(j)}(d^{(j)}))$ and $p_{k}^{(j)}:= \chi_{\frac{1}{2}}(\psi_{k}^{(j)}(p))$ in $F_{k}^{(j)}$.\\
Note that $\|p_{k}^{(j)}-\psi_{k}^{(j)}(p)\| \stackrel{k \to \infty}{\longrightarrow} 0$, hence $\pi_{C}((p_{0}^{(j)},p_{1}^{(j)},\ldots))=\bar{\psi}^{(j)}(p)$, where $\pi_{C}:\prod F_{k}\to C$ denotes the quotient map.

\noindent
Next, we show that there is $K\in \N$ such that $p_{k}^{(j)} \preceq q_{k}^{(j)}$ for all $j=0, \ldots,n$ and $k\ge K$:\\
Suppose this was not the case. Then, for some $j\in \{0,\ldots,n\}$, there would be a sequence $\tau_{k}$ of tracial states on $F_{k}^{(j)}$ such that 
\[
0<\tau_{r_{k}}(q_{r_{k}}^{(j)})<\tau_{r_{k}}(p_{r_{k}}^{(j)})
\]
for a subsequence $(r_{k})_{k\in \N}$ of $(k)_{k\in \N}$ (from now on, we omit the index $r_{k}$ and write $k$ instead).\\
But 
\[
q_{k}^{(j)} \ge \psi_{k}^{(j)}(d^{(j)})-\frac{\alpha}{2}\, ,
\]
so we obtain
\[
\frac{\alpha}{2} \ge \underline{\lim}_{k} (\tau_{k}(q_{k}^{(j)})-\tau_{k}(p_{k}^{(j)})) +\frac{\alpha}{2}\ge \underline{\lim}_{k}(\tau_{k}(\psi_{k}^{(j)}(d^{(j)}))-\tau_{k}(\psi_{k}^{(j)}(p))) \, .
\]
Then there is a free ultrafilter $\omega \in \beta\N\setminus \N$ such that 
\[
\lim_{\omega} \tau_{k} \psi_{k}^{(j)}(d^{(j)})- \lim_{\omega} \tau_{k} \psi_{k}^{(j)}(p)\le \frac{\alpha}{2}\, .
\]
On the other hand, by \ref{qd-traces} and since the $\psi_{k}^{(j)}$ are approximately multiplicative, the map $\tau^{(j)}: A \to \C$, given by
\[
x \mapsto \lim_{\omega} \tau_{k}\psi_{k}^{(j)}(x)
\]
is a tracial state on $A$. (To see that $\tau^{(j)}$ is nonzero, note that $p_{k}^{(j)}\neq 0\; \forall \, k$, hence $F_{k}^{(j)} \neq 0$ and $\tau_{k}(\be_{F_{k}^{(j)}})=1 \; \forall \, k$. Moreover, $\|\psi_{k}^{(j)}(\be_{A})-\be_{F_{k}^{(j)}}\| \to 0$, so $\tau^{(j)}(\be_{A})=1$.)\\
We now see that $\tau^{(j)}(d^{(j)})-\tau^{(j)}(p)\le \frac{\alpha}{2}$, a contradiction to the definition of $\alpha$, so indeed there is $K\in \N$ such that $p_{k}^{(j)} \preceq q_{k}^{(j)}$ for all $j=0, \ldots,n$ and $k\ge K$.

\noindent
Let $s_{k}^{(j)} \in F_{k}^{(j)}$ be partial isometries such that $s_{k}^{(j)}(s_{k}^{(j)})^{*}= p_{k}^{(j)}$ and $(s_{k}^{(j)})^{*}s_{k}^{(j)}\le q_{k}^{(j)}$. These yield partial isometries $s^{(j)}\in C^{(j)}$ such that 
\[
s^{(j)}(s^{(j)})^{*}= \bar{\psi}^{(j)}(\iota(p))
\]
and
\[
(s^{(j)})^{*} s^{(j)} \le \chi_{\frac{\alpha}{2}}(\bar{\psi}^{(j)}(\iota(d^{(j)})))\, .
\]
We may now define 
\[
v^{(j)}:=\bar{\varphi}^{(j)}(\be_{C^{(j)}})^{\halb} \bar{\sigma}^{(j)}(s^{(j)})
\]
and $v:=\sum_{j=0}^{n}v^{(j)}$; these clearly are well-defined elements in $B$. Now we have
\begin{eqnarray*}
(v^{(j)})^{*}v^{(j)} & = & \bar{\varphi}^{(j)}(\be_{C^{(j)}}) \bar{\sigma}^{(j)}((s^{(j)})^{*}s^{(j)}) \\
& \le & \bar{\varphi}^{(j)}(\be_{C^{(j)}}) \bar{\sigma}^{(j)}(g_{0,\frac{\alpha}{2}}(\bar{\psi}^{(j)}(d^{(j)}))) \\
& = & \bar{\varphi}^{(j)}\bar{\psi}^{(j)}(g_{0,\frac{\alpha}{2}}(d^{(j)}))\\
& \le & \bar{\varphi} \bar{\psi}(g_{0,\frac{\alpha}{2}}(d^{(j)}))\\
& = & \iota(g_{0,\frac{\alpha}{2}}(d^{(j)})) \in (\iota(d^{(j)})B\iota(d^{(j)}))^{-} \subset \overline{\iota(d)B \iota(d)} \, 
\end{eqnarray*}
(see \ref{functions} for the definition of $g_{\alpha,\beta}$), from which follows that 
\[
v^{*}v \in \overline{\iota(d)B\iota(d)}
\]
and that
\[
v^{(j)}(v^{(j')})^{*}=0 \mbox{ if } j \neq j' \, .
\]
As a consequence,
\begin{eqnarray*}
v v^{*}&=& \sum_{i,j=0}^{n} v^{(i)}(v^{(j)})^{*}\\
&=& \sum_{j=0}^{n} v^{(j)}(v^{(j)})^{*}\\
&=& \sum_{j=0}^{n} \bar{\varphi}^{(j)}(\be_{C^{(j)}}) \bar{\sigma}^{(j)}(s^{(j)}(s^{(j)})^{*})\\
&=& \sum_{j=0}^{n} \bar{\varphi}^{(j)}\bar{\psi}^{(j)}(\iota(p))\\
&=& \bar{\varphi}\bar{\psi}(\iota(p))\\
&=& \iota(p) \, .
\end{eqnarray*}
Finally, choose a positive lift $h=(h_{0},h_{1},\ldots)$ of $v^{*}v$ in $\prod_{\N} \overline{dAd}$; let $p' \in \prod_{\N} A$ be the lift of $\iota(p)$ given by the constant sequence $(p,p,\ldots)$. By \cite{Lo}, Theorem 8.2.1, there is a lift $t=(t_{0},t_{1},\ldots)$ of $v$ in $\prod_{\N} A$ such that $t t^{*}\le p'$ and $t^{*}t\le h$. Since $t$ lifts $v$, there is $m \in \N$ such that $\|t_{m}t_{m}^{*}-p\|<1$. As a consequence, $t_{m}t_{m}^{*}$ is invertible in $pAp$ and we may define $u:=(t_{m}t_{m}^{*})^{-\halb}t_{m}$. Then $u$ is a partial isometry in $A$ satisfying $uu^{*}=p$ and $u^{*}u\in \overline{dAd}$.
\end{nproof}

\altbn{\label{perforation}}
\begin{ntheorem}
Let $A$ be a separable, simple and unital $C^{*}$-algebra with real rank zero and decomposition rank $n$ for some finite $n$. Then $A$ satisfies Blackadar's second fundamental comparability property, and so do all matrix algebras $M_{k}(A)$ over $A$. In particular, $K_{0}(A)$ is weakly unperforated. 
\end{ntheorem}

\begin{nproof}
Given $p,q \in A$ with $\tau(p)<\tau(q) \; \forall \, \tau \in T(A)$ and set $\varepsilon:= \min \{\tau(q)-\tau(p) \, | \, \tau \in T(A)\}$. Since $(\tau \mapsto \tau(q)-\tau(p))$ is a nonzero positive continuous function on the compact space $T(A)$, we see that $\varepsilon >0$ (and, clearly, $\varepsilon \le 1$). \\
Apply Corollary \ref{orthogonal-subprojections} (with $\eta=1$, say) to obtain $k\in \N$ and $\gamma>0$ which satisfy the assertion of \ref{orthogonal-subprojections}. Apply Corollary \ref{orthogonal-subprojections} oncemore, this time with $\varepsilon \cdot \gamma$ as $\eta$ to obtain $k' \in \N$ and $\gamma'>0$.\\
But this means there are pairwise orthogonal projections $p_{1},\ldots,p_{2^{k}} \in pAp$ with $p=p_{1} \oplus \ldots \oplus p_{2^{k}}$ and $\tau(p_{i})<\varepsilon \cdot \gamma$ $\forall \, \tau \in T(A)$ $i=1,\ldots,2^{k}$. \\
Furthermore, we obtain pairwise orthogonal projections $q_{0}^{(1)}, \ldots, q_{n}^{(1)} \in qAq$ (recall that $2^{k'}\ge n+1$) such that $\tau(q_{j}^{(1)}) > \gamma \tau(q)$ $\forall \, \tau \in T(A)$, $j=0, \ldots,n$. In particular, $\tau(q_{j}^{(1)})>\varepsilon \cdot \gamma> \tau(p_{1})$ for all $j$, so by Proposition \ref{n+1-domination} there is a projection $p_{1}' \in qAq$ with $p_{1}  \sim p_{1}' \le q_{0}^{(1)} \oplus \ldots \oplus q_{n}^{(1)} \le q $.\\
Next, apply the assertion of \ref{orthogonal-subprojections} to $q-p_{1}'$; this yields mutually orthogonal projections $q_{0}^{(2)}, \ldots, q_{n}^{(2)} \le q- p_{1}'$ such that 
\[
\tau(q_{j}^{(2)})> \gamma \cdot \tau(q-p_{1}')\ge \gamma \cdot \tau(q-p) > \gamma \cdot \varepsilon \, .
\]
Again Proposition \ref{n+1-domination} tells us that there is a projection $p_{2}' $ such that $p_{2} \sim p_{2}' \le q-p_{1}'$. \\
We may proceed inductively; after $2^{k}$ steps we have constructed pairwise orthogonal projections $p_{1}', \ldots, p_{2^{k}}' \le q$ with $p_{i} \sim p_{i}'$, from which $p \le q$ follows. \\
The statement also holds for $M_{k}(A)$, since $M_{k}(A)$ is again simple and has decomposition rank $n$ and real rank zero (cf.\ \cite{BP} and \cite{KW}).
\end{nproof}
\alten

\section{Tracial orthogonalization}

\noindent
Below we provide a  method of `almost orthogonalizing order zero maps w.r.t.\ traces'. This will be useful to deal with $C^{*}$-algebras whose tracial state spaces have compact and zero-dimensional extreme boundary.

\altbn{\label{orthogonalization-1}}
\begin{nprop}
Let $a,b$ be positive normalized elements in a  $C^{*}$-algebra $A$ and  $\beta>0$. Then there is $\theta_{b}(a) \in A$ with $0 \le \theta_{b}(a) \le a $ such that the following holds:\\
(i) $\|b \, \theta_{b}(a)\|< \beta^{\frac{1}{8}}$\\
(ii) if $\tau \in T(A)$ satisfies $\tau(a) \tau(b) < \beta$, then $\tau(a) - \tau(\theta_{b}(a))< \beta^{\frac{1}{4}}$.
\end{nprop}

\begin{nproof}
Define 
\[
\theta_{b}(a):= a - a^{\halb} g_{{\beta^{\frac{1}{4}}}/{2},\beta^{\frac{1}{4}}}(a^{\halb}b a^{\halb})a^{\halb} 
\]
(cf.\ \ref{functions}), then $0 \le \theta_{b}(a) \le a$. We now have 
\begin{eqnarray*}
\|\theta_{b}(a) b\|^{2} & = & \| \theta_{b}(a) b^{2} \theta_{b}(a)\| \\
& \le & \| \theta_{b}(a) b \theta_{b}(a)\| \\ 
& \le & \| \theta_{b}(a) b a^{\halb}\| \\
& \le & \|a^{\halb} b a^{\halb} - g_{{\beta^{\frac{1}{4}}}/2 ,\beta^{\frac{1}{4}}}(a^{\halb}b a^{\halb})a^{\halb} b a^{\halb} \| \\
& \le & \frac{\beta^{\frac{1}{4}}}{2} \\
& < & \beta^{\frac{1}{4}} \, ,
\end{eqnarray*}
Furthermore, if $\tau \in T(A)$ satisfies  $\tau(a) \tau(b) < \beta$, then 
\begin{eqnarray*}
\tau(a^{\halb}g_{{\beta^{\frac{1}{4}}}/{2},\beta^{\frac{1}{4}}}(a^{\halb}b a^{\halb}) a^{\halb}) & \le & \tau(g_{{\beta^{\frac{1}{4}}}/{2},\beta^{\frac{1}{4}}}(a^{\halb}b a^{\halb}))\\
& \le & \frac{1}{\beta^{\frac{1}{4}}} \tau(a^{\halb}b a^{\halb}) \\
& < & \frac{1}{\beta^{\frac{1}{4}}} \min\{\tau(a), \tau(b) \} \\
& \le & \frac{1}{\beta^{\frac{1}{4}}} (\tau(a) \tau(b))^{\halb} \\
& < & \beta^{\frac{1}{4}}\, . 
\end{eqnarray*}
\end{nproof}
\alten

\altbn{\label{tracial-orthogonalization}}
\begin{nlemma}
For $F=M_{r_{1}}\oplus \ldots \oplus M_{r_{s}}$, $\mu>0$ and $\eta>0$ there is $\delta>0$ such that the following holds: \\
Let $A$ be a $C^{*}$-algebra and $\varphi:F \to A$ c.p.\ with $\|\varphi\|\le \mu$ and $\ord \varphi|_{M_{r_{i}}}=0$ for $i=1, \ldots, s$. Then there is a c.p.\ map $\varrho:F \to A$ such that \\
(i) $\ord \varrho|_{M_{r_{i}}} =0 \; \forall \, i$ and $\|\varrho(\be_{M_{r_{i}}})\varrho(\be_{M_{r_{j}}})\| < \eta$ whenever $i \neq j$,\\
(ii) $\varrho(x) \le \varphi(x) \; \forall \, x \in F_{+}$, in particular, $\varrho$ is subordinate to $\varphi$,\\
(iii) if, for some $\tau \in T(A)$, the estimates $\tau(\varphi(\be_{M_{r_{i}}}))\tau(\varphi(\be_{M_{r_{j}}}))<\delta$ hold for all $i\neq j$ for which $\varphi(\be_{M_{r_{i}}})\varphi(\be_{M_{r_{j}}}) \neq 0$, then $\tau(\varphi(\be_{M_{r_{i}}}))- \tau(\varrho(\be_{M_{r_{i}}}))<\eta$, $i=1,\ldots,s$. 
\end{nlemma}

\begin{nproof}
By rescaling $\varphi$ and $\eta$, if necessary, without loss of generality  we may assume that $\mu=1$. If we can prove the lemma for $F=M_{r_{1}}\oplus M_{r_{2}}$, then an iterated application will yield the general case. (To make this iteration work, it is crucial that, at each step, the map $\varrho$ satisfies condition (ii). Thus, if $\varphi$ preserves orthogonality for certain elements, so does $\varrho$.) Therefore, let $F=M_{r_{1}}\oplus M_{r_{2}}$ and $\eta>0$ be given and  suppose $\varphi: F \to A$ is c.p.c.\ with $\ord \varphi|_{M_{r_{i}}}=0$ for $i=1,2$.

\noindent
First, choose $0< \beta < 1$ such that $r_{1}^{2} \cdot r_{2} \cdot \beta^{{\frac{1}{16}}} \le \eta$ and set $\delta:= r_{1} \cdot r_{2} \cdot \beta$; fix pairwise orthogonal minimal projections $e_{1}, \ldots , e_{r_{1}} \in M_{r_{1}}$ and $f_{1}, \ldots, f_{r_{2}} \in M_{r_{2}}$. \\
Set $a_{1}^{(0)}:= \varphi(e_{1})$ and, inductively, define $a_{1}^{(k)}:= \theta_{\varphi(f_{k})}(a_{1}^{(k-1)})$ for $k \in \{1, \ldots, r_{2}\}$ by employing Proposition \ref{orthogonalization-1}. From \ref{induced-order0} we obtain a c.p.c.\ map $\varphi^{(1)}: M_{r_{1}} \to A$ with  strict order zero, $\varphi^{(1)}(e_{1})= a_{1}^{(r_{2})}$ and $\varphi^{(1)}(x) \le \varphi(x)$ for all positive $x \in M_{r_{1}}$; in particular, $\varphi^{(1)}$ is subordinate to $\varphi|_{M_{r_{1}}}$.\\
Next, define $a_{2}^{(0)}:= \varphi^{(1)}(e_{2})$. Again, we define $a_{2}^{(k)}:= \theta_{\varphi(f_{k})}(a_{2}^{(k-1)})$ for $k=1, \ldots, r_{2}$ using Proposition \ref{orthogonalization-1}. As above, from \ref{induced-order0} we obtain a c.p.c.\ order zero map $\varphi^{(2)}: M_{r_{1}} \to A$ with $\varphi^{(2)}(e_{2})=a_{2}^{(r_{2})}$. \\
Proceed inductively to construct c.p.c.\ order zero maps $\varphi^{(l)}: M_{r_{1}} \to A$, $l=1, \ldots, r_{1}$ satisfying $\varphi^{(r_{1})}(x) \le \ldots \le \varphi^{(1)}(x) \le \varphi(x) \; \forall \, x \in (M_{r_{1}})_{+}$. If $\tau \in T(A)$ satisfies $\tau(\varphi(\be_{M_{r_{1}}})) \tau(\varphi(\be_{M_{r_{2}}}))< \delta$ then, for any $l=1, \ldots, r_{1}$ and $k=1, \ldots,r_{2}$, $\tau(a_{l}^{(k)}) \tau (\varphi(f_{k}))< \frac{\delta}{r_{1} \cdot r_{2}} = \beta$. Therefore, by Proposition \ref{orthogonalization-1}, $\tau(a_{l}^{(k-1)})- \tau(a_{l}^{(k)})< \beta^{\frac{1}{4}}$. Furthermore, $\tau(a_{l}^{(r_{2})})- \tau(a_{l+1}^{(1)})< \beta^{\frac{1}{4}}$, so  $\tau(a_{1}^{(0)})- \tau(a_{r_{1}}^{(r_{2})})  < r_{1 } \cdot r_{2} \cdot \beta^{\frac{1}{4}}$ and $\tau(\varphi(\be_{M_{r_{1}}})) - \tau(\varphi^{(r_{1})}(\be_{M_{r_{1}}})) < r_{1}^{2} \cdot r_{2} \cdot \beta^{\frac{1}{4}} \le \eta $, since $\tau(\varphi(\be_{M_{r_{1}}}))= r_{1} \cdot \tau(a_{1}^{(0)})$ and $\tau(\varphi^{(r_{1})}(\be_{M_{r_{1}}})) = r_{1} \cdot \tau(a_{r_{1}}^{(r_{2})})$.\\
Next, observe that $\varphi^{(r_{1})}(e_{l}) \le \varphi^{(l)}(e_{l})=a_{l}^{(r_{2})} \le a_{l}^{(k)}$ for $l=1, \ldots, r_{1}$ and $k=1, \ldots, r_{2}$. This yields the estimate 
\begin{eqnarray*}
\|\varphi^{(r_{1})}(\be_{M_{r_{1}}}) \varphi(\be_{M_{r_{2}}})\| & \le & \sum_{l=1}^{r_{1}} \sum_{k=1}^{r_{2}} \|\varphi^{(r_{1})}(e_{l}) \varphi(f_{k})\| \\
& \le & \sum_{l=1}^{r_{1}} \sum_{k=1}^{r_{2}} \| a^{(k)}_{l} \varphi(f_{k})\|^{\halb} \\
& \stackrel{\ref{orthogonalization-1}}{<} & r_{1} \cdot r_{2} \cdot (\beta^{\frac{1}{8}})^{\halb} \\
& = &  r_{1} \cdot r_{2} \cdot \beta^{\frac{1}{16}} \\
& \le & \eta \, .
\end{eqnarray*}
Define a c.p.\ map $\varrho: F \to A$ by $\varrho|_{M_{r_{1}}}:= \varphi^{(r_{1})}$ and $\varrho|_{M_{r_{2}}}:= \varphi|_{M_{r_{2}}}$, then we are done.
\end{nproof}
\alten

\section{Real rank zero and affine functions on $T(A)$}

\noindent
If $A$ is in a class of $C^{*}$-algebras for which  the Elliott conjecture is verified, it must be possible to read off from the invariant wether or not $A$ has real rank zero. Therefore it makes sense to try to rephrase the condition $\rr A =0$ in terms of the Elliott invariant. \\
In the case where $A$ has finite decomposition rank and $T(A)$ has compact and zero-dimensional extreme boundary, it will follow from the results of this section that $\rr A =0$ if and only if $\iota(K_{0}(A)_{+})$ is dense in $\Aff(T(A))_{+}$, the positive continuous affine functions on $T(A)$. The essential step is Theorem \ref{perforation2}, which says that if $\iota(K_{0}(A)_{+})$ is dense in $\Aff(T(A))_{+}$, then the projections of $A$ can be compared in terms of tracial states.

\altbn{\label{halving-projections}}
\begin{nprop}
Let $A$ be a simple unital $C^{*}$-algebra with real rank zero and $\dr A=n<\infty$. For any $\delta >0$ and any projection $p\in A$ there are orthogonal  projections $q, \, q'\in pAp$ which are Murray--von Neumann equivalent and satisfy $\tau(p-q-q')<\delta$ for all $\tau \in T(A)$.
\end{nprop}

\begin{nproof}
First, choose $0<\gamma<1$ and $k \in \N$ (with $2^{k-1}\ge n+1$) such that the assertion of Corollary \ref{orthogonal-subprojections} (with $\eta=1$) holds.\\
Next, choose $m\in \N$ such that $\beta:= \frac{1}{2^{m}}<\gamma$; set $\alpha:= \frac{1}{(9(n+1)^{2})^{m}}$. \\
From a repeated application of Lemma \ref{interpolation} we obtain a projection $p'\le p_{0}:=p$ with 
\[
\alpha \tau(p_{0})<\tau(p')<\beta \tau(p_{0}) < \gamma \tau(p_{0}) \; \forall \, \tau \in T(A)\, . 
\]
By Corollary \ref{orthogonal-subprojections} there are pairwise orthogonal projections $e_{1},\ldots,e_{2(n+1)}\in pAp$ with $\gamma \tau(p)<\tau(e_{i})$ for all $i=1, \ldots, 2(n+1)$ and all $\tau \in T(A)$. By Proposition \ref{n+1-domination}, there are projections $q_{1}\le e_{1}\oplus \ldots \oplus e_{n+1}$ and $q'_{1}\le e_{n+2}\oplus \ldots \oplus e_{2(n+1)}$ such that $p'\sim q_{1}$ and $p' \sim q'_{1}$.

\noindent
Replacing $p$ by $p_{1}:=p-(q_{1}+q'_{1})$ (but not changing $\alpha$, $\beta$ and $\gamma$), we may repeat this construction to obtain orthogonal projections $q_{2}$ and $q'_{2}$ in $p_{1}Ap_{1}$ such that $q_{2}\sim q'_{2}$ and 
\[
\alpha \tau(p_{1})<\tau(q_{2})=\tau(q'_{2})<\beta \tau(p_{1})
\]
for all $\tau \in T(A)$. Induction now yields sequences of projections $(p_{l})_{\N}$, $(q_{l})_{\N}$ and $(q'_{l})_{\N}$ in $pAp$ satisfying:
\begin{itemize}
\item distinct projections in $\{q_{1},q'_{1},q_{2},q'_{2},\ldots \}$ are orthogonal
\item $p_{l+1}=p_{l}-(q_{l+1}+q'_{l+1}) \; \forall \, l$
\item $q_{l} \sim q'_{l} \; \forall \, l$
\item $\alpha \tau(p_{l}) < \tau(q_{l+1})=\tau(q'_{l+1})<\beta \tau(p_{l}) \; \forall \, l \in \N, \, \tau \in T(A)$.
\end{itemize}
Using that, for any $\tau \in T(A)$,
\begin{eqnarray*}¥
\sum_{l=1}^{j+1} \tau(q_{l} +q'_{l}) & > & 2 \alpha \tau(p_{j}) + \sum_{l=1}^{j} \tau(q_{l} + q'_{l}) \\
& = & 2 \alpha \left(\tau(p_{0})- \sum_{l=1}^{j-1} \tau(q_{l} + q'_{l}) \right) + \sum_{l=1}^{j} \tau(q_{l} + q'_{l}) \\
& = & 2 \alpha \tau(p_{0}) + (1-2\alpha) \sum_{l=1}^{j-1} \tau(q_{l}+q'_{l}) +\tau(q_{j}+q'_{j}) \\
& > & 2 \alpha \tau(p_{0}) + (1-2 \alpha) \sum_{l=1}^{j} \tau(q_{l}+q'_{l})
\end{eqnarray*}
one checks  that
\begin{eqnarray*}
\sum_{l=1}^{j} \tau(q_{l}+q'_{l}) & > & 2 \alpha \cdot \tau(p_{0})+ (1-2 \alpha) \cdot 2 \alpha \cdot \tau(p_{0}) + (1-2\alpha)^{2}\cdot 2 \alpha \cdot \tau(p_{0})\\
& & + \ldots + (1-2\alpha)^{j-1} \cdot 2 \alpha \cdot \tau(p_{0}) \\
&=& \left( \sum_{l=0}^{j-1}(1-2\alpha)^{l} \right) \cdot 2 \alpha \cdot \tau(p)\\
&\stackrel{j\to \infty}{\longrightarrow} & \left(\frac{1}{1-(1-2\alpha)}\right) \cdot 2 \alpha \cdot \tau(p)\\
&=&\tau(p)\, .
\end{eqnarray*}
We may regard the projections $h_{j}:= \sum_{l=1}^{j} q_{l}+q'_{l}$ and $p$ as positive continuous functions on the compact metrizable space $T(A)$. By the preceding argument, the $h_{j}$ converge to $p$ pointwise, and since the $h_{j }$ are increasing, we see from Dini's theorem that $h_{j} \to p$ uniformly (as functions on $T(A)$). Thus there is $\bar{\jmath} \in \N$ such that $\tau(p-h_{\bar{\jmath}})<\delta \; \forall \, \tau \in T(A)$. Set $q:=\sum_{l=1}^{\bar{\jmath}}q_{l}$ and $q':=\sum_{l=1}^{\bar{\jmath}}q'_{l}$, then $h_{\bar{\jmath}}= q+q'$ and $q \sim q'$ by our construction.
\end{nproof}
\alten

\altbn{\label{rr0>affine-functions}}
\begin{ncor}
Let $A$ be a simple unital $C^{*}$-algebra with real rank zero and $\dr A=n<\infty$. Then, for any $\varepsilon>0$ and any positive function $f\in \Aff(T(A))$, there is  $[p] \in K_{0}(A)_{+}$ such that $|f(\tau) - \tau(p)|<\varepsilon$. In other words, the image of $K_{0}(A)_{+}$  in $\Aff(T(A))_{+}$ under the natural map $\iota$ is uniformly dense. \\
If, moreover,  $\partial_{e}T(A)$ is compact, then the restriction of $\iota(K_{0}(A)_{+})$ to $\partial_{e}(T(A))$ is dense in $\Ch(\partial_{e}(T(A)))_{+}$.
\end{ncor}

\begin{nproof}
First, define 
\[
N:=\{\iota(x)/2^{k}\, | \, x \in K_{0}(A),\, k\in \N \} \, .
\]
Now \cite{Go}, Theorem 7.9 says that $N$ is uniformly dense in $\Aff(T(A))$ (recall from \ref{r_A} that $S(K_{0}(A))  \approx T(A)$). But, as $A$ satisfies Blackadar's fundamental comparability property (cf.\ Theorem \ref{perforation}), it is easy to check that in fact $N_{+}:=\{\iota(x)/2^{k}\, | \, x \in K_{0}(A)_{+},\, k\in \N \}$ is dense in $\Aff(T(A))_{+}$. On the other hand, $\iota(K_{0}(A)_{+})$ is dense in $N_{+}$ by Proposition \ref{halving-projections}.\\
The second statement holds, since by \cite{Go}, Corollary 7.5,  restriction to $\partial_{e}T(A)$ yields a surjective map: $\Aff(T(A))\twoheadrightarrow \Ch(\partial_{e}(T(A)))$. 
\end{nproof}
\alten

\altbn{\label{orthogonal-functions}}
\begin{nprop}
Let $X$ be a compact and metrizable space with $\dim X =0$. Let $0 \le f_{0}, \ldots, f_{n} \in \Ch(X)$ and $\eta>0$ be given.\\
Then there are $0\le g_{0}, \ldots, g_{n} \in \Ch(X)$ such that $g_{i} \perp g_{j}$ if $i \neq j$, $g_{j}\le f_{j}$ and 
\[
\frac{1}{n+1} \cdot \sum_{l=0}^{n} f_{l} - \eta \le \sum_{l=0}^{n} g_{l} \le \frac{1}{n+1} \cdot \sum_{l=0}^{n} f_{l}
\]
for $i, j=0, \ldots,n$. 
\end{nprop}

\begin{nproof}
Set $h:= \frac{1}{n+1} \cdot \sum_{l=0}^{n} f_{l} -\eta$ and define $U_{j}:=\{t \in X \, | \, f_{j}(t) > h(t)\} $. Obviously, $U_{0}, \ldots, U_{n}$ is an open covering of $X$; since $\dim X=0$, there is an open covering $V_{0},\ldots,V_{n}$ of $X$ such that $V_{j}\subset U_{j}$ and $V_{i} \cap V_{j}= \emptyset$ for all $i\neq j \in \{0,\ldots, n\}$. Now define continuous functions $g_{j}\in \Ch(X)_{+}$ by 
\[
g_{j}(t):= \left\{ 
\begin{array}{ll}
\max \{h(t),\, 0\} & \mbox{if }t\in V_{j}\\
0 & \mbox{else} \, .
\end{array}
\right.
\]
It is straightforward to check that these functions have the right properties.
\end{nproof}
\alten

\altbn{\label{dividing-elements}}
\begin{nprop}
Let $A$ be a separable, simple, unital $C^{*}$-algebra with $\dr A=n<\infty$ and such that $\partial_{e} T(A)$ is compact and satisfies $\dim \partial_{e} T(A)=0$. Furthermore, suppose that $\iota(K_{0}(A)_{+}) \subset \Ch(\partial_{e} T(A))_{+}$ is dense. \\
Given a positive normalized element $d \in A$ and $\gamma>0$, there are $n+1$ pairwise orthogonal normalized elements $0 \le d_{0}, \ldots, d_{n} \in \overline{dAd}$ such that $\frac{1}{(n+1)^{2}} \tau(d)-\gamma<\tau(d_{i})$ for all $i=0, \ldots, n$ and $\tau \in T(A)$.    
\end{nprop}

\begin{nproof}
Choose a system $(F_{k},\psi_{k},\varphi_{k})_{\N}$ of c.p.\ approximations for $A$ such that the $\psi_{k}$ are approximately multiplicative and the $\varphi_{k}$ are $n$-decomposable with respect to $F_{k}=\bigoplus_{j=0}^{n}F_{k}^{(j)}$; denote the respective decompositions of $\psi_{k}$ and $\varphi_{k}$ by $\psi_{k}^{(j)}$ and $\varphi_{k}^{(j)}$, $k\in \N$, $j=0,\ldots,n$. \\
By Proposition \ref{simple-nf-approximations} we may even assume that the irreducible representations of $h_{k}F_{k}h_{k}$, where $h_{k}:= \chi_{\frac{\gamma}{16}}(\psi_{k}(d))\in F_{k}$, are at least $k$-dimensional. \\
Set $C:=M_{n+1}\otimes M_{n+1}$ and choose $*$-homomorphisms $\pi_{k}:C\to h_{k}F_{k}h_{k}$ which are maximal in the sense that, for each $k$, the algebra $(h_{k}-\pi_{k}(\be_{C}))F_{k}(h_{k}-\pi_{k}(\be_{C}))$ does not contain a copy of $C$. \\
If $\tau$ is a tracial state on $A$, then the maps $\tau \verk \varphi_{k}$ are positive contractive traces on $F_{k}$ (because the $\varphi_{k}$ are $n$-decomposable; cf.\ \ref{order0-traces}), and so are the restrictions $\tau \verk \varphi_{k}|_{h_{k}F_{k}h_{k}}$. But then, since the $\pi_{k}$ are maximal and the irreducible representations of $h_{k}F_{k}h_{k}$ are at least $k$-dimensional, it is clear that $\tau\verk \varphi_{k}(h_{k}-\pi_{k}(\be_{C})) < \frac{(n+1)^{2}}{k}$. As a consequence, there is $K_{0}\in \N$ such that 
\begin{eqnarray}
\tau \verk \varphi_{k}\verk \pi_{k}(\be_{C}) &\ge& \tau \verk \varphi_{k}(h_{k}) - \frac {\gamma}{32} \nonumber\\
&\ge& \tau \verk \varphi_{k} (f_{\frac{\gamma}{16},\frac{\gamma}{8}}(\psi_{k}(d)))-\frac {\gamma}{32} \nonumber\\
&\ge& \tau \verk \varphi_{k}(\psi_{k}(d))-\frac {\gamma}{16}-\frac {\gamma}{32} \label{fd-trace-estimate}\\
&\ge& \tau  (d) -\frac{\gamma}{16}- 2\frac {\gamma}{32} \nonumber\\
&=&\tau(d)-\frac{\gamma}{8} \nonumber
\end{eqnarray}
for all $\tau \in T(A)$ and $k\ge K_{0}$ (see \ref{functions} for the definition of $f_{\alpha,\beta}$).  

\noindent
With $\frac{\gamma}{8 \cdot k}$ in place of $\eta$, $\mu=n+1$ and $C\otimes \C^{n+1}$ in place of $F$ we apply Lemma \ref{tracial-orthogonalization} to obtain $\delta_{k}>0$; we may clearly assume  $\delta_{k}<\frac{\gamma}{16(n+1)}$.\\
Choose a rank-one projection $e\in M_{n+1}$ and, for each $k\in \N$ and $j=0,\ldots,n$, define 
\[
f_{k}^{(j)}:= \varphi_{k}^{(j)}\verk \pi_{k}^{(j)}(\be_{M_{n+1}}\otimes e) \in A \, .
\]
Now, by our assumption on the space of exremal tracial states on $A$, from Proposition \ref{orthogonal-functions} (with $\frac{\gamma}{8}$ in place of $\eta$ and $(f_{k}^{(j)}- \frac{\gamma}{16}\cdot \be)_{+}$ in place of $f_{j}$) we obtain positive functions $g_{k}^{(j)} \in \Ch(\partial_{e}T(A))$ satisfying
\begin{itemize}
\item $g_{k}^{(i)} \perp g_{k}^{(j)}$ if $i\neq j$
\item $g_{k}^{(j)}(\tau)\le \frac{1}{n+1} \tau((f_{k}^{(j)}-\frac{\gamma}{16})_{+})$
\item $\frac{1}{n+1}\sum_{l=0}^{n}\tau((f_{k}^{(l)}-\frac{\gamma}{16})_{+})-\frac{\gamma}{8}\le \sum_{l=0}^{n} g_{k}^{(l)}(\tau)\le \frac{1}{n+1}\sum_{l=0}^{n}\tau((f_{k}^{(l)}-\frac{\gamma}{16})_{+})$
\end{itemize}
for all $i,j \in \{0,\ldots,n\}$ and $\tau\in T(A)$.\\
By our assumption, the image of $K_{0}(A)_{+}$ in $\Ch(\partial_{e}T(A))_{+}$ is dense, so for each $k$ there is $m_{k}\in \N$ and a projection $q_{k}^{(j)} \in M_{m_{k}}(A)$ such that 
\[
g_{k}^{(j)}(\tau)- \frac{\delta_{k}}{3(n+1)^{2}}<\tau(q_{k}^{(j)})<g_{k}^{(j)}(\tau)+ \frac{\delta_{k}}{3(n+1)^{2}}
\]
for all $k\in \N$, $j=0,\ldots,n$ and $\tau \in T(A)$. \\
Let $e_{0},\ldots,e_{n} \in M_{n+1}$ be a set of pairwise orthogonal rank-one projections, then
\[
f_{k}^{(j)}= \sum_{l=0}^{n}\varphi_{k}^{(j)}\verk \pi_{k}^{(j)}(e_{l}\otimes e)
\]
and
\[
\tau \verk \varphi_{k}^{(j)}\verk \pi_{k}^{(j)}(e_{l}\otimes e)= \frac{1}{n+1}\tau(f_{k}^{(j)})\, .
\]
Now 
\[
\tau(q_{k}^{(j)}) \le \tau \verk \varphi_{k}^{(j)} \verk \pi_{k}^{(j)}(e_{l}\otimes e)
\]
for all $l=0,\ldots,n$. Since $\ord(\varphi_{k}^{(j)})=0$, we see from Proposition \ref{n+1-domination} that each $q_{k}^{(j)}$ is Murray--von Neumann equivalent to a projection in $(f_{k}^{(j)}Af_{k}^{(j)})^{-}$. Misusing our notation we may assume that $q_{k}^{(j)} \in (f_{k}^{(j)}Af_{k}^{(j)})^{-}$. \\
From \ref{induced-order0}, for each $k\in \N$ and $j\in \{0,\ldots, n\}$ we obtain a c.p.c.\ map $\bar{\varphi}_{k}^{(j)}:M_{n+1} \to A$ with strict order zero (and the same supporting $*$-homomorphism as $\varphi_{k}^{(j)} \verk \pi_{k}^{(j)}$) such that $\bar{\varphi}_{k}^{(j)}(e)=q_{k}^{(j)}$. \\
The $\bar{\varphi}_{k}^{(j)}$ add up to c.p.\ maps $\bar{\varphi}_{k}:M_{n+1}\otimes \C^{n+1} \to A$ with $\|\bar{\varphi}_{k}\|\le n+1$. Next, we may apply Lemma \ref{tracial-orthogonalization} to find c.p.\ maps $\varrho_{k}:\bigoplus_{j=0}^{n} M_{n+1}\to A$ (we denote the components by $\varrho_{k}^{(j)}$) satisfying \\
(i) $\ord \varrho^{(j)}_{k}=0$ and $\|\varrho^{(i)}_{k}(\be_{M_{n+1}}) \varrho^{(i)}_{k}(\be_{M_{n+1}})\| < \frac{\gamma}{8 \cdot k}$ for all $i \neq j$\\
(ii) $\varrho_{k}(x)\le \bar{\varphi}_{k}(x) \; \forall \, x \in (M_{n+1}\otimes \C^{n+1})_{+}$\\
(iii) if, for some $\tau \in T(A)$, the estimates $\tau(\bar{\varphi}_{k}^{(i)}(\be_{M_{n+1}})) \tau(\bar{\varphi}_{k}^{(j)}(\be_{M_{n+1}}))< \delta_{k} $ hold whenever $i\neq j$, then $\tau(\bar{\varphi}_{k}^{(j)}(\be_{M_{n+1}}))-\tau(\varrho_{k}^{(j)}(\be_{M_{n+1}}))<\frac{\gamma}{8 \cdot k}$.

\noindent
But now we check that, if $i\neq j$, then
\begin{eqnarray*}
\lefteqn{\tau(\bar{\varphi}_{k}^{(i)}(\be_{M_{n+1}}))\tau(\bar{\varphi}_{k}^{(j)}(\be_{M_{n+1}}))}\\
& = & (n+1)^{2} \tau(q_{k}^{(i)})\tau(q_{k}^{(j)})\\
&\le&(n+1)^{2} \left(g_{k}^{(i)}(\tau)+\frac{\delta_{k}}{3(n+1)^{2}} \right) \left(g_{k}^{(j)}(\tau)+\frac{\delta_{k}}{3(n+1)^{2}} \right)\\
&\le& (n+1)^{2} \cdot 3 \, \frac{\delta_{k}}{3(n+1)^{2}} \\
&<&\delta_{k} \, ,
\end{eqnarray*}
so indeed
\[
0\le \tau(\bar{\varphi}_{k}^{(j)}(\be_{M_{n+1}})- \varrho_{k}^{(j)}(\be_{M_{n+1}}))<\frac{\gamma}{8 \cdot k}
\]
for all $k \in \N$, $j=0,\ldots,n$ and all $\tau\in T(A)$. 

\noindent
By (\ref{fd-trace-estimate}) there is $K_{0}\in \N$ such that, for all $k\ge K_{0}$ and $\tau \in T(A)$,
\begin{eqnarray*}
 \tau(d) & \le & \frac{1}{(n+1)^{2}} \left(\tau \verk \varphi_{k} \verk\pi_{k}(\be_{C}) + \frac{\gamma}{8} \right) \\
&=& (n+1)\sum_{l=0}^{n}\tau(f_{k}^{(l)}) + \frac{\gamma}{8}  \\
&\le& (n+1)\sum_{l=0}^{n} \tau((f_{k}^{(l)}-\frac{\gamma}{16})_{+}) +(n+1)\frac{\gamma}{16} +\frac{\gamma}{8} \\
&\le&  (n+1)^{2} \sum_{l=0}^{n} g_{k}^{(l)}(\tau) + (n+1)^{2} \frac{\gamma}{8} + (n+1)\frac{\gamma}{16} +\frac{\gamma}{8} \\
& \le &  (n+1)^{2} \sum_{l=0}^{n} \tau(q_{k}^{(l)}) +(n+1)^{3}\frac{\delta}{3(n+1)^{2}} + (n+1)^{2} \frac{\gamma}{8} + (n+1)\frac{\gamma}{16} +\frac{\gamma}{8} \\ 
& \le &  (n+1)^{2}\sum_{l=0}^{n} \tau(\varrho_{k}^{(l)}(e_{0})) +(n+1)^{2}\frac{\gamma}{8(n+1)}  +(n+1)^{3}\frac{\delta}{3(n+1)^{2}} \\
& & + (n+1)^{2} \frac{\gamma}{8} + (n+1)\frac{\gamma}{16} +\frac{\gamma}{8}  \, .
\end{eqnarray*}
As a consequence, we see that
\begin{eqnarray*}
\lefteqn{\tau(\varrho_{k} \verk \kappa(e_{0}))}\\
& \ge & \frac{1}{(n+1)^{2}} \tau(d) - \left( \frac{\gamma}{8(n+1)}+\frac{\delta}{3(n+1)}+\frac{\gamma}{8}+ \frac{\gamma}{16(n+1)}+\frac{\gamma}{8(n+1)^{2}} \right)\\
&\ge & \frac{1}{(n+1)^{2}}\tau(d) - \frac{\gamma}{2} \, ,
\end{eqnarray*}
where $\kappa : M_{n+1} \to \bigoplus_{0,\ldots ,n}  M_{n+1}$ denotes the diagonal embedding. It follows from our construction that 
\[
\|\varrho_{k}(x) g_{0,\frac{\gamma}{16}}(d) - \varrho_{k}(x)\| \to 0 
\]
and that 
\[
\|\varrho_{k}^{(i)}(\be_{M_{n+1}}) \varrho_{k}^{(j)}(\be_{M_{n+1}})\| \to 0
\]
for $i \neq j$, whence the $\varrho_{k}$ induce an order zero map 
\[
\varrho: \bigoplus_{0, \ldots, n} M_{n+1} \to \prod_{\N} \overline{dAd} /\bigoplus_{\N} \overline{dAd}
\]
(again, see \ref{functions} for the definition of $g_{\alpha,\beta}$). \\
Since order zero maps are semiprojective (cf.\ \cite{Wi2}),  there are $K_{1}\in \N$ (we may assume $K_{1}\ge K_{0}$) and an order zero map 
\[
\bar{\varrho}: \bigoplus_{0, \ldots, n} M_{n+1} \to \prod_{k\ge K_{1}} \overline{dAd}
\]
with 
\[
\|\bar{\varrho}_{k} \verk \kappa (e_{0})-\varrho_{k} \verk \kappa (e_{0})\|< \frac{\gamma}{2} \; \forall \, k \ge K_{1} \, .
\]
Then in particular we have $\bar{\varrho}_{K_{1}}\verk \kappa (M_{n+1}) \in \overline{dAd}$ and, for $j=0, \ldots,n$ and $\tau \in T(A)$, from \ref{order0-traces} we obtain 
\[
\tau(\bar{\varrho}_{K_{1}} \verk \kappa (e_{j}))= \tau(\bar{\varrho}_{K_{1}} \verk \kappa (e_{0})) \ge \frac{1}{(n+1)^{2}} \tau(d) - \gamma \, .
\]
The $\bar{\varrho}_{K_{1}} \verk \kappa (e_{0}), \ldots,\bar{\varrho}_{K_{1}} \verk \kappa (e_{n})$ are pairwise orthogonal, since $\bar{\varrho}_{K_{1}} \verk \kappa$ has strict order zero.
\end{nproof}
\alten

\altbn{\label{perforation2}}
\begin{ntheorem}
Let $A$ be a separable, simple, unital $C^{*}$-algebra with $\dr A=n<\infty$ and such that $\partial_{e} T(A)$ is compact and satisfies $\dim \partial_{e} T(A)=0$. Furthermore, suppose that $\iota(K_{0}(A)_{+}) \subset \Ch(\partial_{e} T(A))_{+}$ is dense. \\
If $p\in A$ is a projection and $d\in A$ is a positive normalized element such that $\tau(p)<\tau(d)$ for all $\tau$ in $T(A)$, then $p$ is Murray--von Neumann equivalent to a projection in $\overline{dAd}$. \\
In particular, $A$ satisfies Blackadar's second fundamental comparability property and $K_{0}(A)$ is weakly unperforated.
\end{ntheorem}

\begin{nproof}
We may clearly assume that $p$ is nonzero, for otherwise there is nothing to show.  Set 
\[
\alpha:= \min\{\tau(d)-\tau(p)\, |\, \tau \in T(A)\} \, .
\] 
Since $T(A)$ is compact and $\tau(p)<\tau(d)$ for all $\tau$, and because $d$ is nonzero and normalized, we see that $0<\alpha < 1$. 

\noindent
Next, choose $k\in \N$ such that $1-\frac{k}{k+2(n+1)^{2}}<\frac{\alpha}{2(n+1)^{2}}$. With  $\beta:=\frac{1}{k+2(n+1)^{2}}$ we have 
\[
0<1-k\cdot \beta< \frac{\alpha}{2(n+1)^{2}} \mbox{ and } \beta<\frac{1-(k-1)\cdot \beta}{2(n+1)^{2}}\, .
\]
We may regard $\beta\cdot p$ as a continuous function on $\partial_{e}T(A)$ (defined by $\tau \mapsto \beta \cdot \tau(p)$). Since the image of $K_{0}(A)_{+}$ in $\Ch(\partial_{e}T(A))_{+}$ is dense, it is straightforward to see that there exists a projection $e \in M_{r}(A)$ (for some $r\in \N$) such that 
\[
0<\tau(p)-k \tau(e)<\frac{\alpha}{2(n+1)^{2}}\tau(p)
\]
and
\[
\tau(e)<\frac{\tau(p)-(k-1)\tau(e)}{2(n+1)^{2}}\; \forall \, \tau \in T(A)\, .
\]
Because $p\neq 0$, the number 
\[
\gamma:= \min\{\frac{\tau(p)}{2(n+1)^{2}}\, | \, \tau \in T(A)\} 
\]
is nonzero and  by Proposition \ref{dividing-elements} there are pairwise orthogonal positive normalized elements $d_{0}, \ldots,d_{n} \in pAp$ satisfying 
\[
\tau(d_{j})\ge \frac{\tau(p)}{(n+1)^{2}}-\gamma \ge \frac{\tau(p)}{2(n+1)^{2}}>\tau(e)
\] 
for all $j=0,\ldots, n$ and $\tau \in T(A)$. By Proposition \ref{n+1-domination} there is a projection $p_{0}\in pAp$ which is Murray--von Neumann equivalent to $e$. \\
In the same manner (but replacing $p$ by $p-p_{0}$) we obtain from Propositions \ref{dividing-elements} and  \ref{n+1-domination} a projection $p_{1} \le p-p_{0}$ which is Murray--von Neumann equivalent to $e$. \\
Using that 
\[
\tau(e)< \frac{\tau(p-(p_{0}+\ldots +p_{i}))}{2(n+1)^{2}}
\] 
for all $i\le k-1$, we can proceed inductively to construct pairwise orthogonal projections $p_{0}, \ldots, p_{k-1} \in pAp$ which are all Murray--von Neumann equivalent to $e$. Set $p_{k}:=p-(p_{0}+\ldots+p_{k-1})$, then 
\[
\tau(p_{i})<\frac{\alpha}{2(n+1)^{2}} \tau(p)
\]
for all $i=0,\ldots, k$ and $\tau \in T(A)$.

\noindent
Now choose an approximate unit $0\le h_{0}\le h_{1}\le \ldots$ for $\overline{dAd}$ satisfying $h_{i}\ge (d-\frac{\alpha}{2})_{+}$ and $h_{i+1}h_{i}=h_{i}$ for all $i\in \N$.\\
Since $\tau(h_{0})\ge \frac{\alpha}{2}$ and $\tau(p_{0})\le \frac{\alpha}{2(n+1)^{2}}$ for all $\tau \in T(A)$, it follows from Propositions \ref{dividing-elements} and  \ref{n+1-domination}, that there is a projection $q_{0} \in \overline{h_{0}Ah_{0}}$ with $p_{0}\sim q_{0}$. \\
We now have $0\le h_{1}-q_{0}$ and $\tau(h_{1}-q_{0})\ge \frac{\alpha}{2}$, so we may apply \ref{dividing-elements} and  \ref{n+1-domination} oncemore to obtain a projection $q_{1} \in \overline{(h_{1}-q_{0})A(h_{1}-q_{0})}$ with $p_{1} \sim q_{1}$. Induction yields projections $q_{j}\in \overline{(h_{j}-q_{j-1})A(h_{j}-q_{j-1})}$ with $q_{j} \sim p_{j}$ for $j=1, \ldots, k$.\\
It is clear from our construction that the $q_{j}$ are pairwise orthogonal and sum up to a projection $q_{0}\oplus \ldots \oplus q_{k} \in \overline{dAd}$ which is Murray--von Neumann equivalent to $p$.
\end{nproof}
\alten

\section{Tracial approximations}

\noindent
In this section we recall the notion of tracial rank zero, prove our main result and derive its corollaries.

\altbn{\label{d-tr}}
Recall from \cite{Li1} that a simple unital $C^{*}$-algebra $A$ is said to have tracial rank zero, if the following holds: \\
Given $\varepsilon>0$, $0 \neq a \in A_{+}$ and  $G \subset A$ finite, there is a finite-dimensional $C^{*}$-algebra $B \subset A$ such that

\noindent
(i) $\|[\be_{B},x]\|< \varepsilon \; \forall \, x \in B$,\\
(ii) $\dist(pxp,B) < \varepsilon \; \forall \, x \in B$,\\
(iii) $\be_{A}-\be_{B}$ is Murray--von Neumann equivalent to a projection in $\overline{aAa}$.\\

\noindent
It was shown in \cite{Li1} that simple, separable, unital $C^{*}$-algebras with tracial rank zero have real rank zero, stable rank one, are quasidiagonal and  their $K_{0}$-groups are weakly unperforated and have the Riesz interpolation property. In \cite{Li3}, Lin showed that the class $\mathfrak{A}$ of all simple, separable, unital, nuclear $C^{*}$-algebras with tracial rank zero which satisfy the Universal coefficient theorem (cf.\ \cite{Bl1}) indeed satisfies Conjecture \ref{elliott-conjecture}.  
\alten

\altbn{\label{dr-tr}}
\begin{ntheorem}
Let $A$ be a separable, simple, unital $C^{*}$-algebra with $\dr A=n<\infty$. Suppose that $\partial_{e}T(A)$ is compact,  that $\dim \partial_{e}T(A)=0$ and that $\iota(K_{0}(A)_{+}) \subset \Ch(\partial_{e}T(A))_{+}$ is dense. Then $A$ has tracial rank zero.
\end{ntheorem}

\begin{nproof}
Let $a_1, \ldots, a_m \in A$ be positive normalized elements and let $\varepsilon>0$ be given. Choose $\eta>0$ such that $\eta^{\frac{1}{4}}<\frac{1}{10(n+1)}$ and such that $2 \eta +4 \eta^{\frac{1}{4}}<\varepsilon$. Next, choose a c.p.\ approximation $(F=M_{r_{1}}\oplus \ldots \oplus M_{r_{s}}, \psi, \varphi)$ for $a_{1}, \ldots, a_{m},a_{1}^{2},\ldots,a_{m}^{2}$ within $\eta$ such that $\varphi$ is $n$-decomposable with respect to $\{1,\ldots,s\}= \coprod_{j=0}^{n} I_{j}$ and such that $\|\psi(a_{l})^{2}-\psi(a_{l}^{2})\|<\eta$ for all $l$. \\
Denote the restrictions $\varphi|_{M_{r_{i}}}$ by $\varphi_{i}$ and let $\sigma_{i}:M_{r_{i}} \to A''$ be supporting $*$-homo\-morphisms for $\varphi_{i}$. For each $i$ let $\{e_{k,l}^{(i)}\, |\, k,l=1,\ldots, r_{i}\}$ be a set of matrix units for $M_{r_{i}}$. Define $B_{i}:=(\sigma_{i}(e_{1,1}^{(i)})A\sigma_{i}(e_{1,1}^{(i)}))^{-}$.\\
Following \ref{induced-order0}, we then obtain $*$-homomorphisms $\bar{\sigma}_{i}:B_{i}\otimes M_{r_{i}} \to  A''$ and c.p.c.\ maps $\bar{\varphi}_{i}:B_{i}\otimes M_{r_{i}} \to A$ given by
\[
\bar{\sigma}_{i}(a\otimes x):= \sum_{k=1}^{r_{i}} \sigma_{i}(e_{k,1}^{(i)})a \sigma_{i}(e_{1,k}^{(i)}) \sigma_{i}(x) 
\]
and
\[
\bar{\varphi}_{i}(a \otimes x):= \varphi_{i}(\be_{M_{r_{i}}})^{\halb} \bar{\sigma}_{i}(a \otimes x)\varphi_{i}(\be_{M_{r_{i}}})^{\halb} \, .
\]
In particular we see that the $\bar{\varphi}_{i}$ add up to a c.p.c.\ map 
\[
\bar{\varphi}: \bigoplus_{i=1}^{s} B_{i}\otimes M_{r_{i}} \to A
\]
satisfying $\bar{\varphi}\verk \kappa=\varphi$, where $\kappa:F \to \bigoplus_{i=1}^{s}B_{i}\otimes M_{r_{i}}$ is the obvious unital embedding.\\
For each $i$ define a c.p.c.\ map $\varphi'_{i}:M_{r_{i}} \to A$ by $\varphi'_{i}(x):= f_{\eta^{\frac{1}{4}},2\eta^{\frac{1}{4}}}(\varphi_{i}(\be_{M_{r_{i}}}))\sigma_{i}(x) $ (see \ref{functions} for the definition of $f_{\alpha,\beta}$). The $\varphi_{i}' $ add up to a c.p.c.\ map $\varphi':F \to A$ which is subordinate to $\varphi$ and satisfies $\|\varphi'(x)- \varphi(x)\|\le (n+1) \eta^{\frac{1}{4}}\|x\| $ for all $x \in F_{+}$. \\
For $j=0, \ldots,n$ define $h^{(j)}:= \sum_{i\in I_{j}} \varphi'(\be_{M_{r_{i}}})$. Regarding the $h^{(j)}$ as continuous functions on $\partial_{e}T(A)$, we may apply Proposition \ref{orthogonal-functions} to obtain positive functions $g^{(0)}\oplus \ldots \oplus g^{(n)} \in \Ch(\partial_{e}T(A))$ such that $g^{(j)} \perp g^{(j')}$ if $j \neq j'$, $g^{(j)}(\tau)\le \tau(h^{(j)})$ and 
\[
\frac{1}{n+1}\sum_{j=0}^{n} \tau(h^{(j)}) - \frac{1}{10(n+1)} \le \sum_{j=0}^{n} g^{(j)}(\tau) \le \frac{1}{n+1} \sum_{j=0}^{n} \tau(h^{(j)})
\]
for all $\tau \in \partial_{e}T(A)$. For $i=1, \ldots, s$ define $g_{i}:= \inf \{g^{(j)},\, \iota(\varphi'(\be_{M_{r_{i}}}))\} \in \Ch(\partial_{e}T(A))$, where $j$ is such that $i\in I_{j}$.

\noindent
Since the $g^{(j)}$ are pairwise orthogonal, we have for all  $\tau \in \partial_{e}T(A)$ that
\[
\sum_{j=0}^{n} g^{(j)}(\tau) \le \sum_{i=1}^{s} g_{i}(\tau) \, .
\]
We now  obtain
\begin{eqnarray*}
\frac{1}{n+1} \tau(\varphi'(\be_{F})) - \frac{1}{10(n+1)} & = & \frac{1}{n+1} \sum_{j=0}^{n} \tau(h^{(j)}) - \frac{1}{10(n+1)}\\
& \le & \sum_{j=0}^{n} g^{(j)}(\tau) \\
& \le & \sum_{i=1}^{s} g_{i}(\tau)\\
& \le & \tau(\varphi'(\be_{F}))\, .
\end{eqnarray*}
Furthermore, we have $g_{i}(\tau) \le \tau(\varphi'(\be_{M_{r_{i}}})) \; \forall \, i$ and $\sum_{i\in I_{j}} g_{i} \perp \sum_{i\in I_{j'}} g_{i}$ whenever $j \neq j' \in \{0,\ldots,n\}$.

\noindent
For $\gamma:=\frac{\eta}{s}$, there is $\alpha>0$ as in  Proposition \ref{weakly-stable}.  
Obtain $0<\delta< \gamma$ from Lemma \ref{tracial-orthogonalization} (with $\mu = n+1$ and $\eta':= \min\{ \eta, \alpha \}$ in place of $\eta$). Since $\iota(K_{0}(A)_{+})$ is dense in $\Ch(\partial_{e}T(A))$, there are $r\in \N$ and projections $p_{i}\in M_{r}(A)$ such that
\[
\frac{1}{r_{i}}(g_{i}(\tau)- \delta/2)< \tau(p_{i})<\frac{1}{r_{i}}(g_{i}(\tau) +\delta/2)\; \forall \, \tau \in \partial_{e}T(A)\, .
\]
But then $\tau(p_{i})<\tau(\varphi'(e_{1,1}^{(i)}))$ for all $\tau$ and $i$, so by Theorem \ref{perforation2} there are projections $p_{i}' \in (\varphi_{i}'(e_{1,1}^{(i)})A\varphi_{i}'(e_{1,1}^{(i)}) )^{-}$ with $p_{i}' \sim p_{i} \; \forall \, i$. \\
From \ref{induced-order0} we obtain order zero maps $\varphi_{i}'':M_{r_{i}} \to A$ with $\varphi_{i}''(e_{1,1}^{(i)})=p_{i}'$ and such that $\sigma_{i}$ is a supporting $*$-homomorphism for $\varphi_{i}'' \; \forall \, i$ (in particular this means that the sum $\varphi''$ is subordinate to  $\varphi$). \\
Now if $i\neq i' \in \{1,\ldots,s\}$ are so that $\varphi_{i}''(\be_{M_{r_{i}}})\varphi_{i'}''(\be_{M_{r_{i'}}}) \neq 0$, then 
\[
\varphi_{i}(\be_{M_{r_{i}}})\varphi_{i'}(\be_{M_{r_{i'}}})\neq 0 \, ,
\]
hence $i \in I_{j}$ and $i' \in I_{j'}$ for distinct $j$ and $j'$. But then $g_{i} \perp g_{i'}$, so we obtain
\[
\tau(\varphi_{i}''(\be_{M_{r_{i}}}))\tau(\varphi_{i'}''(\be_{M_{r_{i'}}})) <   \delta
\]
for all $\tau \in \partial_{e}T(A)$. \\
From Lemma \ref{tracial-orthogonalization} we now obtain a c.p.\ map $\varrho: F \to A$ which has strict order zero on the matrix blocks of $F$, satisfies $\|\varrho(\be_{M_{r_{i}}}) \varrho(\be_{M_{r_{i'}}})\|< \eta' \le \alpha$ for $i \neq i'$ and which is subordinate to $\varphi''$ (hence to $\varphi$). 
Furthermore, $\varrho(x) \le \varphi''(x) \; \forall \, x \in F_{+}$ and $\tau(\varphi''(\be_{F}))- \tau(\varrho(\be_{F}))< \eta' \; \forall \, \tau \in \partial_{e}T(A)$. By Proposition \ref{weakly-stable} there is an order zero map $\bar{\varrho}: F \to A$ with $\|\bar{\varrho}(x) - \varrho(x)\| < \gamma \|x\| \; \forall \, x \in F$.\\
As above, we find projections $q_{i}\in (\bar{\varrho}(e_{1,1}^{(i)})A \bar{\varrho}(e_{1,1}^{(i)}))^{-}$ with 
\[
\frac{1}{r_{i}}(\tau(\bar{\varrho}(\be_{M_{r_{i}}}))-\gamma)<\tau(q_{i})<\frac{1}{r_{i}}(\tau(\bar{\varrho}(\be_{M_{r_{i}}}))+ \gamma) \, ,
\]
from which follows that 
\[
\frac{1}{r_{i}}(\tau({\varrho}(\be_{M_{r_{i}}}))-2 \gamma)<\tau(q_{i})<\frac{1}{r_{i}}(\tau({\varrho}(\be_{M_{r_{i}}}))+ 2 \gamma) \, .
\]
This in turn yields a c.p.c.\ map $\varrho': F \to A$ with the  properties that $\varrho'$ is subordinate to $\bar{\varrho}$ (in particular, it has strict order zero and the restrictions $\varrho'_{i}:=\varrho'|_{M_{r_{i}}}$ have supporting $*$-homomorphisms $\sigma_{i}$) and 
\[
\varrho'(e_{1,1}^{(i)})= q_{i} \in (\bar{\varrho}(e_{1,1}^{(i)})A \bar{\varrho}(e_{1,1}^{(i)}))^{-} \, .
\]
Since $\varrho'$ has strict order zero and maps projections to projections, by \cite{Wi1}, Proposition 3.2.b), it is a $*$-homomorphism . We proceed to check that the finite-dimensional subalgebra $\varrho'(F) \subset A$ has the right properties. First, note that for all $\tau \in \partial_{e}T(A)$
\begin{eqnarray*}
\tau(\varrho'(\be_{F}))& = & \sum_{i=1}^{s} r_{i} \cdot \tau(q_{i}) \\
& \ge & \sum_{i=1}^{s} \tau(\varrho(\be_{M_{r_{i}}})) - 2 \cdot s \cdot \gamma \\
& \ge & \tau(\varphi''(\be_{F})) - 2 \cdot s \cdot \gamma - \eta'\\
& > & \sum_{i=1}^{s} g_{i}(\tau) - 2 \cdot s \cdot \gamma - \eta' - s \cdot \delta\\
& \ge & \frac{1}{n+1} \tau(\varphi'(\be_{F})) - \frac{1}{10(n+1)} - s \cdot (2 \cdot \gamma + \delta) - \eta'\\
& \ge & \frac{1}{n+1} \tau(\varphi(\be_{F})) - \eta^{\frac{1}{4}}-  \frac{1}{10(n+1)} - s \cdot (2 \cdot \gamma + \delta) - \eta'\\
& \ge & \frac{1}{n+1} \tau(\be) - \frac{\eta}{n+1} - \eta^{\frac{1}{4}}-  \frac{1}{10(n+1)} - s \cdot (2 \cdot \gamma + \delta) - \eta'\\
& \ge & \frac{1}{n+1} - \frac{\eta}{n+1} - \eta^{\frac{1}{4}} - \frac{1}{10(n+1)} - 4 \eta \\
& > & \frac{1}{n+1} - \frac{1}{2(n+1)}\\
& = & \frac{1}{2(n+1)} \, .
\end{eqnarray*}
Using functional calculus it is straightforward to construct positive elements $c_{i}\in C^{*}(\varphi_{i}(e_{1,1}^{(i)}))$ satisfying $c_{i}\varphi_{i}(e_{1,1}^{(i)})\varphi_{i}'(e_{1,1}^{(i)})= \varphi_{i}'(e_{1,1}^{(i)})$ and $\|c_{i}\|\le \frac{1}{\eta^{\frac{1}{4}}}$; we may regard the $c_{i}$ as elements of $B_{i}$  and set $c:=c_{1} \oplus \ldots \oplus c_{s} \in \bigoplus B_{i}$. \\
Recall that $\varrho'$ is subordinate to $\varphi'$, so in particular we have $q_{i} = \varrho'_{i}(e_{1,1}^{(i)}) \in (\varphi'_{i}(e_{1,1}^{(i)})A\varphi'_{i}(e_{1,1}^{(i)}))^{-}$, from which follows that $c_{i}^{\halb}\varphi_{i}(e_{1,1}^{(i)})^{\halb}q_{i}=q_{i}$. Now for an arbitrary element $x = x_{1} \oplus \ldots \oplus x_{s} \in F$ we have
\begin{eqnarray*}
\bar{\varphi}(\kappa(x)c) & = & \bar{\varphi}(c\kappa(x)) \\
& = & \bar{\varphi}(c^{\halb}\kappa(x)c^{\halb})\\
& = & \sum_{i=1}^{s} \bar{\varphi}_{i}((c_{i}^{\halb}q_{i}c_{i}^{\halb})\otimes x_{i})\\
& = & \sum_{i=1}^{s} \varphi_{i}(\be_{F})^{\halb} (\sum_{k=1}^{r_{i}} \sigma_{i}(e_{k,1}^{(i)}) c_{i}^{\halb} q_{i} c_{i}^{\halb} \sigma_{i}(e_{1,k}^{(i)})) \sigma_{i}(x_{i}) \varphi_{i}(\be_{F})^{\halb}\\
& = & \sum_{i=1}^{s} \sum_{k=1}^{r_{i}} (\sigma_{i}(e_{k,1}^{(i)})\varphi_{i}(e_{1,1}^{(i)})^{\halb}c_{i}^{\halb}q_{i}c_{i}^{\halb}\varphi_{i}(e_{1,1}^{(i)})^{\halb}\sigma_{i}(e_{1,k}))\sigma_{i}(x_{i})\\
& = & \sum_{i=1}^{s} \sum_{k=1}^{r_{i}} \sigma_{i}(e_{k,1}^{(i)})q_{i}\sigma_{i}(e_{1,k}^{(i)})\sigma_{i}(x_{i})\\
& = & \sum_{i=1}^{s}\varrho_{i}'(x_{i})\\
& = & \varrho'(x)\, ,
\end{eqnarray*}
which in particular shows that $\bar{\varphi}(c)= \varrho'(\be_{F})$.\\
Using that $\bar{\varphi}\verk \kappa= \varphi$ and that $\psi$ is almost multiplicative on the $a_{l}$ we obtain for $l=1, \ldots, m$
\begin{eqnarray*}
\|\bar{\varphi}(\kappa(\psi(a_{l}))^{2})- \bar{\varphi}(\kappa(\psi(a_{l})))^{2}\| & = & \|\bar{\varphi}(\kappa(\psi(a_{l})^{2}))- \bar{\varphi}(\kappa(\psi(a_{l})))^{2}\| \\
& = & \|\varphi(\psi(a_{l})^{2})- \varphi (\psi(a_{l}))^{2}\| \\
& \le & \|\varphi \psi (a_{l}^{2})- \varphi \psi(a_{l})^{2}\| + \eta \\
& < & 4 \eta \, .
\end{eqnarray*}
As a consequence, from \cite{KW}, Lemma 3.5, we see that for $l=1, \ldots, m$ and any positive $b \in \bigoplus_{i=1}^{s} (B_{i}\otimes M_{r_{i}})$ 
\[
\|\bar{\varphi}(b \kappa(\psi(a_{l})))- \bar{\varphi}(b) \bar{\varphi}(\kappa(\psi(a_{l})))\| \le \sqrt{4 \eta} \|b\| = 2 \eta^{\halb} \|b\| \, .
\]
In particular we obtain 
\begin{eqnarray*}
\|\varrho'(\be_{F})a_{l} - a_{l} \varrho'(\be_{F})\| & \le & \|\bar{\varphi}(c) \varphi \psi(a_{l}) - \varphi \psi(a_{l}) \bar{\varphi}(c) \| + 2 \eta \\
& = & \|\bar{\varphi} (c) \bar{\varphi}(\kappa(\psi(a_{l}))) - \bar{\varphi}(\kappa(\psi(a_{l})))\bar{\varphi}(c)\| + 2 \eta\\
& \le & \|\bar{\varphi}(c \kappa(\psi(a_{l}))) - \bar{\varphi}(\kappa(\psi(a_{l}))c)\| + 2 \eta + 4 \eta^{\halb}\cdot \|c\| \\
& \le & 0 + 2 \eta + 4 \eta^{\halb} \frac{1}{\eta^{\frac{1}{4}}}\\
& < & \varepsilon 
\end{eqnarray*}
and 
\begin{eqnarray*}
\|\varrho'(\be_{F})a_{l} - \varrho'(\psi(a_{l}))\| & \le & \|\bar{\varphi}(c) \bar{\varphi}(c \kappa(\psi(a_{l})))\| + \eta\\
& \le & 2 \eta^{\halb} \cdot \eta^{-\frac{1}{4}} + \eta \\
& < & \varepsilon \, . 
\end{eqnarray*}
\end{nproof}
\alten

\altbn{\label{dichotomy}}
\begin{ncor}
Let $A$ be a separable, simple, unital $C^{*}$-algebra with $\dr A=n<\infty$, $\partial_{e}T(A)$  compact and   $\dim \partial_{e}T(A)=0$. \\
Then $\rr A=0$ iff   $\iota(K_{0}(A)_{+}) \subset \Ch(\partial_{e}T(A))_{+}$ is dense iff $A$ has tracial rank zero. \\
If $A$ has a unique trace $\tau$, then $A$ has tracial rank zero iff $\tau(K_{0}(A)_{+})$ is dense in $\R_{+}$; otherwise, $A$ is stably isomorphic to a simple unital $C^{*}$-algebra containing no nontrivial projections.
\end{ncor}

\begin{nproof}
The first statement follows from Corollary \ref{rr0>affine-functions}, Theorem \ref{dr-tr} and \ref{d-tr}.\\
If $A$ has only one tracial state $\tau$, then either $\tau(K_{0}(A)_{+}) \subset \R_{+}$ is dense (in which case $\rr A =0$ by Theorem \ref{dr-tr} and \ref{d-tr}), or $\tau(K_{0}(A)_{+}) \setminus \{0\}$ has a minimum, say $\tau([p])$, where $p \in A\otimes \Kh$ is  some nonzero projection. But then $p(A\otimes \Kh) p$ does not contain any nontrivial projection ($\tau([p])$ is the minimum), and $(p(A\otimes \Kh)p) \otimes \Kh \cong A\otimes \Kh$ by Brown's theorem.
\end{nproof}
\alten

\altbn{\label{classification}}
\begin{ncor}
For $i=1,2$, let $A_{i}$ be  separable, simple, unital $C^{*}$-algebras with $\dr A_{i}<\infty$, $\rr A_{i} = 0$, $\partial_{e}T(A_{i})$  compact and   $\dim \partial_{e}T(A_{i})=0$; suppose that $A_{1}$ and $A_{2}$ satisfy the Universal coefficient theorem. Then:

\noindent
(i)  $A_{1} \cong A_{2}$ if and only if 
\[
(K_{0}(A_{1}),K_{0}(A_{1})_{+},[\be_{A_{1}}],K_{1}(A_{1})) \cong (K_{0}(A_{2}),K_{0}(A_{2})_{+},[\be_{A_{2}}],K_{1}(A_{2})) \, .
\]
(ii) $A$ and $B$ are $AH$-algebras of topological dimension at most 3.

\noindent
(iii) $A$ and $B$ are $ASH$-algebras of topological dimension at most 2. In particular, $\dr A, \dr B \le 2$.
\end{ncor}

\begin{nproof}
(i) follows from Theorem \ref{dr-tr}, Corollary \ref{dichotomy} and the classification theorem of \cite{Li3} (cf.\ \ref{d-tr}). \\
(ii) is a combination of results from Dadarlat, Elliott and Gong (cf.\ \cite{Ro}, Theorem 3.3.5).\\
(iii) follows from results of Elliott (cf.\ \cite{Ro}, Theorem 3.4.4) and \cite{Wi3}, Example 1.10.
\end{nproof}
\alten

\altbn{\label{tensorZ}}
In \cite{JS}, Jiang and Su constructed a simple, unital, projectionless $C^{*}$-algebra $\Zh$, which is $KK$-equivalent to $\C$. The following Corollary (at least partially) generalizes \cite{JS}, Theorem 5, where it was shown that simple unital $AF$-algebras are invariant under tensoring with $\Zh$.\\ 
\begin{ncor}
Let $A$ be a separable, simple, unital $C^{*}$-algebra with $\dr A=n<\infty$, $\rr A = 0$, $\partial_{e}T(A)$  compact and   $\dim \partial_{e}T(A)=0$; suppose that $A$ satisfies the Universal coefficient theorem. Then $A \cong A \otimes \Zh$.
\end{ncor}

\begin{nproof}
Since $\dr \Zh =1$, we have  $\dr (A \otimes \Zh)< \infty $ by \cite{KW}, 3.3.
$A$ and $A \otimes \Zh$ both satisfy the UCT, so by the K\"unneth Theorem (cf.\ \cite{Bl1}) the canonical embedding $\alpha: A \hookrightarrow A \otimes \Zh$ induces an isomorphism of groups $\alpha_{*}: K_{*}(A) \cong K_{*}(A \otimes \Zh)$. Since $\Zh$ has a unique tracial state,  the same embedding induces a homeomorphism $\bar{\alpha}:T(A) \approx T(A \otimes \Zh)$. Furthermore, $\rr (A)=0 $, hence $\iota(K_{0}(A)_{+}) \subset \Ch(\partial_{e}T(A))_{+}$ is dense.  As a consequence,  $\iota \verk \alpha_{*} (K_{0}(A)_{+}) \subset \iota(K_{0}(A \otimes \Zh)_{+}) \subset \Ch(\partial_{e}T(A\otimes \Zh))_{+}$ is dense and $\rr A \otimes \Zh =0$ by Corollary \ref{dichotomy}. Therefore, the maps $r_{A}$ and $r_{A\otimes \Zh}$ are homeomorphisms (cf.\ \ref{r_A}).   Now by Theorem \ref{perforation2} $K_{0}(A)$ and $K_{0}(A \otimes \Zh)$ satisfy Blackadar's second fundamental comparability property (cf.\ \ref{comparability}), hence 
\begin{eqnarray*}
K_{0}(A)_{+} &= & \{x \in K_{0}(A) \, | \, r_{A}(\tau)(x)>0 \; \forall \, \tau \in T(A)\} \\
& = & \{\alpha_{*}(x)  \in K_{0}(A \otimes \Zh) \, | \, r_{A \otimes \Zh}(\bar{\alpha}(\tau))(\alpha_{*}(x))>0 \; \forall \, \tau \in T(A) \}\\
& = & \{x \in K_{0}(A \otimes \Zh) \, | \, r_{A \otimes \Zh}(\tau(x))>0 \; \forall \, \tau \in T(A \otimes \Zh)\}\\
& = & K_{0}(A \otimes \Zh)_{+} \, . 
\end{eqnarray*}
This shows that in fact $K_{0}(A) \cong K_{0}(A \otimes \Zh)$ as ordered groups, so $A \cong A \otimes \Zh$ by Corollary \ref{classification}.
\end{nproof}
\alten

\altbn{\label{Villadsens-example}}
\begin{nexamples}
(i) There are simple separable $C^{*}$-algebras which are approximately homogeneous (hence stably finite and quasidiagonal) but have infinite decomposition rank:\\  
In \cite{Vi1}, Villadsen has constructed simple, separable, unital, approximately homogeneous $C^{*}$-algebras with infinite topological dimension (as $AH$-algebras) and stable rank strictly larger than $1$. It follows from Villadsen's construction that these algebras have unique tracial states and contain projections which are arbitrarily small in trace, so their $K_{0}$-groups are mapped to dense subsets of $\R$ by the tracial states. Now Corollary \ref{dichotomy} implies that Villadsen's examples have infinite decomposition rank. \\
(ii) Finally, we remark  that the conditions of Theorem \ref{dr-tr} are satisfied for many $C^{*}$-algebras. Let $X$ be a compact metrizable space. Then the space $M_{1}^{+}(X)$ of positive regular Borel probability measures on $X$ is a Choquet simplex with extreme boundary $X$ and there is an isometric isomorphism between $\Ch(X)$ and $\Aff(M_{1}^{+}(X))$ (\cite{Go}, Corollaries 10.18 and 11.20). Take any countable dense subgroup $G_{0}$ of $\Ch(X)$ containing $\be_{X}$; with the strict ordering this becomes a partially ordered simple abelian group which is weakly unperforated and has the Riesz interpolation property.  The inclusion $G_{0} \hookrightarrow \Aff(M_{1}^{+}(X))$ induces a continuous affine map $M_{1}^{+}(X) \approx S(\Aff(M_{1}^{+}(X))) \to S(G_{0})$ which is bijective, since $G_{0}\subset \Aff(M_{1}^{+}(X)))$ is dense. Therefore, $S(G_{0})$ is homeomorphic to $M_{1}^{+}(X)$ and  $\partial_{e}S(G_{0}) \approx X$. Finally, fix a countable abelian group $G_{1}$. By results of Elliott (cf.\ \cite{El2}) in connection with \cite{Wi3}, 1.10, there is a simple unital $ASH$ algebra $A$ with decomposition rank at most 2 which has  $(G_{0},(G_{0})_{+},\be_{X},M_{1}^{+}(X),G_{1})$ as its Elliott invariant. In particular, $\partial_{e} T(A) \approx X$ and $\iota(K_{0}(A)_{+})$ is dense in $\Ch(\partial_{e} T(A))$. If $X$ is zero-dimensional,  $A$ satisfies the conditions of Theorem \ref{dr-tr}.     
\end{nexamples}
\alten





\end{document}